\documentclass[a4paper,fleqn]{cas-dc}

\usepackage[numbers]{natbib}
\usepackage{orcidlink}
\usepackage{algorithm}
\usepackage{algorithmicx}
\usepackage{algpseudocode}
\makeatletter
\newenvironment{breakablealgorithm}
{
		\begin{center}
			\refstepcounter{algorithm}
			\hrule height.8pt depth0pt \kern2pt
			\renewcommand{\caption}[2][\relax]{
				{\raggedright\textbf{\ALG@name~\thealgorithm} ##2\par}%
				\ifx\relax##1\relax 
				\addcontentsline{loa}{algorithm}{\protect\numberline{\thealgorithm}##2}%
				\else 
				\addcontentsline{loa}{algorithm}{\protect\numberline{\thealgorithm}##1}%
				\fi
				\kern2pt\hrule\kern2pt
			}
		}{
		\kern2pt\hrule\relax
	\end{center}
}
\makeatother

\def\tsc#1{\csdef{#1}{\textsc{\lowercase{#1}}\xspace}}
\tsc{WGM}
\tsc{QE}

\begin{document}
\let\WriteBookmarks\relax
\def\floatpagepagefraction{1}
\def\textpagefraction{.001}

\shorttitle{}    

\shortauthors{}  

\title [mode = title]{TRUST-TAEA: A trustworthiness-guided two-archive evolutionary algorithm with variable-grouping sparse search for large-scale multi-objective optimization}  

\tnotemark[1] 

\tnotetext[1]{} 

%

\author[1, 2]{Junyi Cui}[style=chinese]
\author[1,2]{Chao Min}[style=chinese]
\author[1,3]{Stanis{\l}aw Mig{\'o}rski}
\cormark[1]
\author[1,2]{Binrong Wang}[style=chinese]
\author[4]{Yonglan Xie}[style=chinese]
\address[1]{School of Sciences, Southwest Petroleum University, Chengdu, 610500, China}
\address[2]{Institute for Artificial Intelligence, Southwest Petroleum University, Chengdu, 610500, China}
\address[3]{Jagiellonian University in Krakow, Faculty of Mathematics and Computer Science, Krakow, 30348, Poland}
\address[4]{School of Geoscience and Technology, Southwest Petroleum University, Chengdu, 610500, China}

\cortext[1]{Jagiellonian University in Krakow, Faculty of Mathematics and Computer Science, 30348 Krakow, Poland. E-mail address: stanislaw.migorski@uj.edu.pl.} 

\begin{abstract}
Large-scale multi-objective optimization problems (LSMOPs) remain challenging due to the high-dimensional decision spaces, complex variable interactions, and limited function evaluation budgets, which make it difficult to balance the convergence, diversity, and stability. Existing two-archive evolutionary algorithms can alleviate the conflict between convergence and diversity, but they often underuse archive reliability and problem-structure information, leading to inefficient search, incomplete front coverage, and late-stage archive drift. To address these issues, this paper proposes TRUST-TAEA, a trustworthiness-guided two-archive evolutionary algorithm. Archive trustworthiness is defined by integrating evolutionary progress with convergence-archive maturity, and is used to coordinate variable-grouping sparse search, anchor-probing compensatory search, and archive stabilization. TRUST-TAEA is evaluated on the LSMOP benchmark suite with 500--5000 decision variables and 2, 3-objectives. Experimental results show that TRUST-TAEA achieves superior and highly competitive performance in terms of convergence, diversity, and stability. A three-objective day-ahead scheduling case of a grid-connected microgrid further demonstrates its practical applicability, where TRUST-TAEA obtains the best IGD$^+$ value and generates a feasible dispatch strategy balancing cost, emissions, and grid-power fluctuation.
\end{abstract}


\begin{keywords}
 Large-scale\sep Multi-objective optimization\sep Trustworthiness guidance\sep Sparse search\sep Two-archive evolutionary algorithm\sep
\end{keywords}

\maketitle

\section{Introduction}\label{sec1}
Multi-objective optimization plays an important role in engineering design \cite{bib1,bib2,bib3}, resource scheduling \cite{bib4,bib5,bib6}, and intelligent decision-making \cite{bib7,bib8}. Its goal is to obtain a set of Pareto-optimal solutions that provide different trade-offs among multiple conflicting objectives. With the increasing complexity of practical optimization tasks, LSMOPs have become an important research topic in evolutionary multi-objective optimization \cite{bib9,bib10}. Compared with the conventional multi-objective optimization problems, LSMOPs usually involve hundreds or thousands of decision variables, leading to a much larger search space, stronger variable interactions, and more limited effective search under a fixed function-evaluation budget. As a result, it is difficult for an algorithm to maintain convergence, diversity, and robustness simultaneously. In particular, many function evaluations may be consumed by ineffective full-dimensional perturbations, while variables with different roles may contribute unequally to Pareto-front shape and convergence behavior. Therefore, improving search efficiency in high-dimensional decision spaces while maintaining a well-distributed approximation front remains a key challenge in LSMOP research.

Existing studies on LSMOPs can be roughly divided into two categories. The first category focuses on exploiting problem-structure information to reduce search complexity. Representative strategies include decision-variable analysis and variable grouping \cite{bib11,bib12}, cooperative co-evolution \cite{bib13, bib14}, problem transformation or reformulation \cite{bib15, bib16}, and sparse or pattern-mining-based search \cite{bib17, bib18}. These methods improve the probability of effective search by decomposing, transforming, or selectively activating decision variables. The second category improves the evolutionary framework itself, such as adaptive offspring generation \cite{bib19}, learning-assisted search acceleration \cite{bib20}, reference-vector guidance \cite{bib21}, and archive updating strategies \cite{bib22}. Although these studies have improved the scalability of evolutionary algorithms from different perspectives, a common issue remains: the search behavior is often not sufficiently coordinated with the reliability of the current population or archive information. In high-dimensional problems, using immature search guidance too aggressively may cause premature convergence, whereas insufficient exploitation of reliable information may waste evaluations in later stages.

Against this background, two-archive evolutionary algorithms have attracted increasing attention because they maintain solution sets with different functional roles. The two-archive framework in C-TAEA \cite{bib23} uses one archive to emphasize convergence and another archive to preserve diversity and complementary search information. Related two-archive studies, such as Two\_Arch2 \cite{bib24} and two-archive matching-based optimization \cite{bib25}, further demonstrate that separating convergence pressure from diversity maintenance can improve the balance between convergence and distribution. Compared with single-archive methods, the two-archive mechanism provides a more flexible framework for preserving promising convergent solutions while retaining exploratory information in less represented regions.

Recent studies have further extended two-archive algorithms by designing archive updating criteria \cite{bib26}, archive coordination strategies \cite{bib27}, and search mechanisms \cite{bib28}. For example, Hu et al. \cite{bib29} proposed a two-archive model for multimodal multi-objective optimization, Xu et al. \cite{bib30} introduced reinforced knowledge sharing into a two-archive framework, and Chen et al. \cite{bib31} developed a two-population two-archive evolutionary framework for constrained multi-objective optimization. However, most existing two-archive improvements are mainly designed for general multi-objective optimization scenarios. When directly applied to LSMOPs, they still face several limitations.

\begin{figure}
	\centering
	\includegraphics[width=\linewidth]{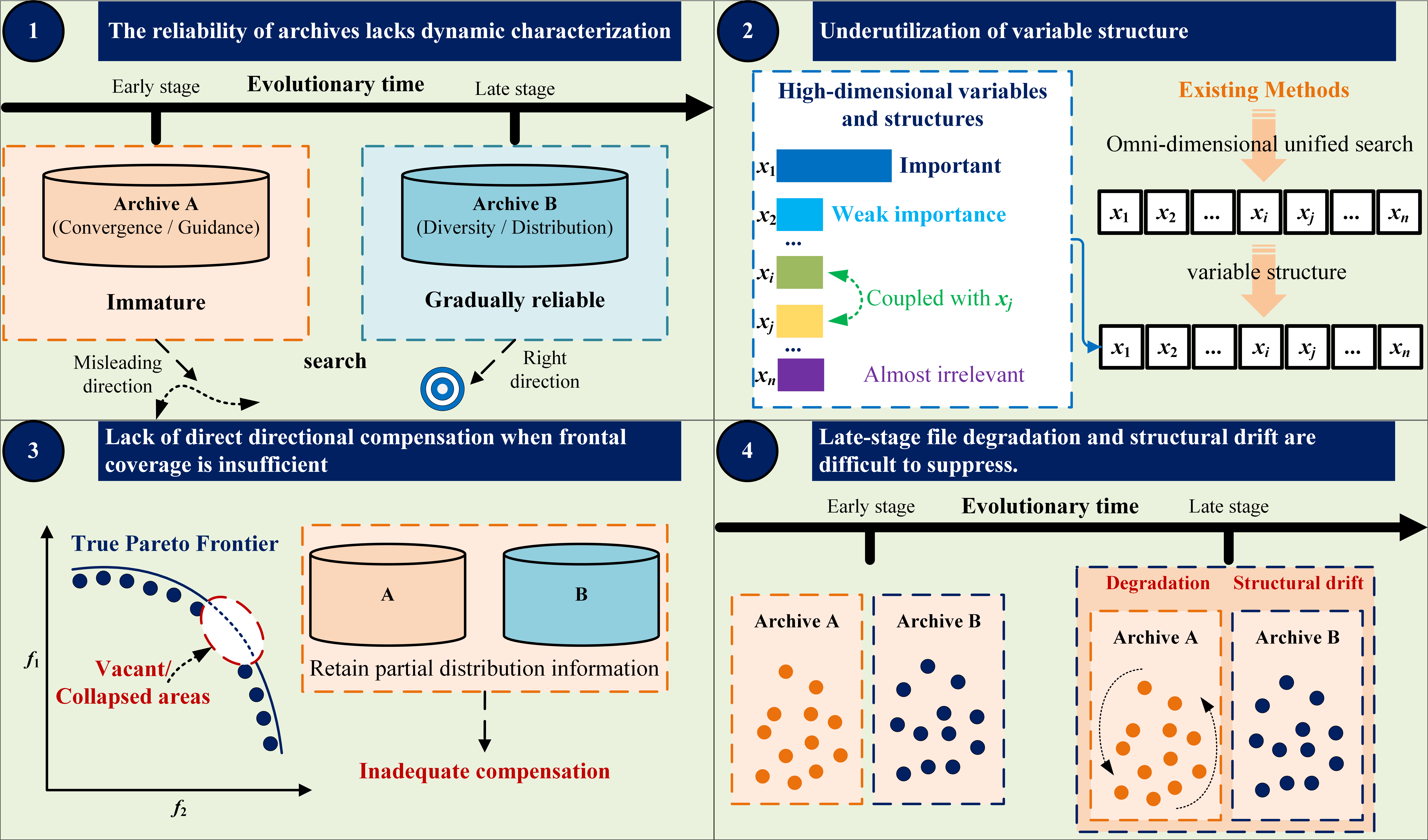} 
	\caption{Shortcomings of existing two-archive improvements in large-scale high-dimensional problems.}
	\label{fig1}
\end{figure}

First, existing two-archive algorithms usually assume that archive information can provide useful guidance throughout the evolutionary process, but they rarely characterize whether the convergence archive is sufficiently reliable at a given stage. In early evolution, an immature archive may contain unstable or poorly distributed solutions, and using it as strong guidance may mislead the search. Second, most two-archive frameworks still rely on full-dimensional reproduction or unified variation operators, which are inefficient for LSMOPs with imbalanced variable contributions and strong structural coupling. Third, although the preference archive can preserve diversity to some extent, it does not explicitly provide targeted compensation for undercovered regions of the Pareto front. Finally, recent LSMOP studies have shown that performance stability is also important in large-scale optimization \cite{bib32}, yet archive degeneration and structural drift in the later evolutionary stages are still not sufficiently addressed in existing two-archive frameworks, as illustrated in Fig.~\ref{fig1}.

Motivated by these considerations, this paper proposes TRUST-TAEA, a trustworthiness-guided two-archive evolutionary algorithm for solving LSMOPs. In TRUST-TAEA, archive trustworthiness is defined by jointly considering evolutionary progress and convergence-archive maturity, and is used as a unified control quantity to regulate offspring generation. Based on this trustworthiness measure, a variable-grouping sparse search mechanism activates only a subset of informative variable groups and performs guided evolution with structural repair, thereby improving the effectiveness of high-dimensional search. Meanwhile, an anchor-probing compensatory search mechanism directionally injects structurally feasible candidate solutions into undercovered front regions to enhance distribution completeness. Furthermore, a checkpoint-based archive stabilization strategy is introduced to suppress archive degeneration and structural drift in later stages. The main contributions of this paper are summarized as follows.
\begin{itemize}
	\item An archive trustworthiness-guided regulation mechanism is proposed to dynamically coordinate exploration and exploitation by jointly considering evolutionary progress and convergence-archive maturity.
	
	\item A variable-grouping sparse search strategy is designed to activate informative variable groups and improve the proportion of effective search in high-dimensional decision spaces.
	
	\item An anchor-probing compensatory search mechanism is developed to provide targeted completion for undercovered regions of the Pareto front.
	
	\item A checkpoint-based archive stabilization strategy is introduced to preserve reliable archive structures and improve late-stage robustness.
\end{itemize}

The remainder of this paper is organized as follows. Section~\ref{sec2} briefly reviews the background of LSMOPs with decision variable decomposition and the conventional two-archive evolutionary framework. Section~\ref{sec3} presents the proposed TRUST-TAEA in detail. Section~\ref{sec4} reports the experimental settings, comparative results, and corresponding analysis. Finally, Section~\ref{sec5} concludes the paper and discusses future research directions.

\section{Preliminaries}\label{sec2}
This section briefly reviews the basic characteristics of LSMOPs and the conventional two-archive mechanism.

\subsection{LSMOP with decision-variable decomposition}
Large-scale multi-objective optimization aims to solve optimization problems with high-dimensional decision spac\-es and multiple conflicting objectives. Its general formulation can be written as
\begin{flalign}\label{1}
	& \min_{\boldsymbol{x}\in\Omega\subset\mathbb{R}^D} F(\boldsymbol{x})
	=
	\left(f_1(\boldsymbol{x}),f_2(\boldsymbol{x}),\ldots,f_M(\boldsymbol{x})\right), &
\end{flalign}
where $\boldsymbol{x}=(x_1,x_2,\ldots,x_D)$ denotes the $D$-dimensional decision vector ($ D> 100$), $\Omega$ is the feasible decision space, and $M$ is the number of objectives. As a result, full-dimensional variation may spend too many function evaluations on the variables that contribute little to the convergence or front distribution. Therefore, the key difficulty is not only the size of the search space, but also the unequal and coupled effects of different decision variables.

A common way to improve the search efficiency is to decompose the decision vector into several variable groups:
\begin{flalign}\label{2}
	& \mathcal{G}=\{G_1,G_2,\ldots,G_K\}, \;
	\boldsymbol{x}=
	\left(
	\boldsymbol{x}^{(1)},
	\boldsymbol{x}^{(2)},
	\ldots,
	\boldsymbol{x}^{(K)}
	\right), &
\end{flalign}
where $G_k\subseteq\{1,2,\ldots,D\}$ denotes the index set of the $k$-th variable group, $\boldsymbol{x}^{(k)}=\{x_j\mid j\in G_k\}$, and $\bigcup_{k=1}^{K}G_k=\{1,2,\ldots,D\}$. In this way, the original high-dimensional problem can be transformed into several lower-dimensional subproblems, and the search effort can be allocated more selectively.

Existing studies usually construct $\mathcal{G}$ from two perspectives. The first is interaction- or contribution-based grouping, where variables with similar effects on the objectives or strong interaction relationships are assigned to the same group \cite{bib33,bib34}. In a general form, the grouping can be described as
\begin{flalign}\label{3}
	& G_k=\left\{j \mid \eta_j \in \mathcal{I}_k \right\}, \; k=1,2,\ldots,K, &
\end{flalign}
where $\eta_j$ denotes the estimated contribution or interaction feature of variable $x_j$, and $\mathcal{I}_k$ is the corresponding feature interval or cluster. This type of method attempts to reduce interference among unrelated variables and improve the effectiveness of local search.

The second is sparsity- or importance-driven grouping, where only a small number of influential groups are activated during each generation \cite{bib35}. Let $\pi_k^t$ denote the activation probability of group $G_k$ at generation $t$. The active group set can be expressed as
\begin{flalign}\label{4}
	& \mathcal{G}_a^t \subseteq \mathcal{G}, \;
	|\mathcal{G}_a^t|=K_a^t, \;
	\Pr(G_k\in\mathcal{G}_a^t)=\pi_k^t. &
\end{flalign}
The corresponding active decision-variable index set is
\begin{flalign}\label{5}
	& \mathcal{I}^t=\bigcup_{G_k\in\mathcal{G}_a^t}G_k, &
\end{flalign}
then the search operator only modifies variables in $\mathcal{I}^t$:
\begin{flalign}\label{6}
	& x_j^{t+1} =
	\begin{cases}
		\mathcal{E}_j(x_j^t), & j\in \mathcal{I}^t,\\
		x_j^t, & j\notin \mathcal{I}^t,
	\end{cases} &
\end{flalign}
where $\mathcal{E}_j(\cdot)$ denotes the variable-level search operator. This sparse search form is especially useful for LSMOPs because it increases the proportion of effective perturbations under limited function-evaluation budgets.

These variable-decomposition ideas provide the structural basis for TRUST-TAEA. Different from general grouping methods that mainly focus on static dimensionality reduction, TRUST-TAEA further combines variable grouping with archive trustworthiness. That is, the number of active groups, the group-selection probability, and the structural repair strength are dynamically regulated according to the reliability of the current convergence archive.

\subsection{Conventional two-archive evolutionary framework}
TAEAs provide an effective framework for balancing convergence and diversity in multi-objective optimization. Following the two-archive idea in C-TAEA \cite{bib23}, two complementary archives are maintained during evolution. One archive is convergence-oriented and acts as the main driving force toward the Pareto front. The other archive is diversity- or preference-oriented and preserves complementary search information, especially regions that are insufficiently covered by the convergence archive. In this paper, these two archives are denoted by the convergence archive $C^t$ and the preference archive $A^t$, respectively.

Let the population at generation $t$ be denoted by $P^t$, and $C^t$ and $A^t$ denote the two archives with archive capacity $N$. In a conventional two-archive framework, offspring are generated from the population and archives:
\begin{flalign}\label{7}
	& \mathcal{Q}^t=\mathcal{R}(P^t,C^t,A^t), &
\end{flalign}
where $\mathcal{R}(\cdot)$ denotes the reproduction operator. In C-TAEA, the two archives cooperate through a restricted mating selection mechanism, where the parents can be selected from different archives according to their evolutionary status. In a general form, this mating process can be expressed as
\begin{flalign}\label{8}
	& \boldsymbol{y}_i^t=
	\mathcal{R}\left(
	\boldsymbol{x}_i^t,
	\boldsymbol{p}_{c,i}^t,
	\boldsymbol{p}_{a,i}^t
	\right), \;
	\boldsymbol{p}_{c,i}^t\in C^t,\;
	\boldsymbol{p}_{a,i}^t\in A^t, &
\end{flalign}
where $\boldsymbol{x}_i^t\in P^t$ is the current parent solution, $\boldsymbol{p}_{c,i}^t$ and $\boldsymbol{p}_{a,i}^t$ are archive-based mating solutions selected from $C^t$ and $A^t$, respectively, and $\boldsymbol{y}_i^t$ is the generated offspring. This formulation reflects the collaborative role of the two archives: $C^t$ provides the convergence pressure, while $A^t$ supplies the complementary distribution information.

After offspring generation, the current population, offspring, and two archives are combined into a union set:
\begin{flalign}\label{9}
	& U^t=P^t\cup \mathcal{Q}^t\cup C^t\cup A^t. &
\end{flalign}

The convergence archive is then updated by a convergence-oriented selection operator:
\begin{flalign}\label{10}
	& C^{t+1}=\mathcal{S}_C(U^t,N), &
\end{flalign}
where $\mathcal{S}_C(\cdot)$ selects at most $N$ solutions with better convergence quality, usually according to nondominated sorting, convergence indicators, or decomposition-based fitness.

Given the updated convergence archive, the preference archive is updated through
\begin{flalign}\label{11}
	& A^{t+1}=\mathcal{S}_A(U^t,C^{t+1},N), &
\end{flalign}
where $\mathcal{S}_A(\cdot)$ selects complementary solutions that improve diversity, preference distribution, or regions insufficiently represented by $C^{t+1}$. The next-generation population is reconstructed from the two archives:
\begin{flalign}\label{12}
	& P^{t+1}=\mathcal{R}_{\mathrm{pop}}(C^{t+1},A^{t+1},N), &
\end{flalign}
where $\mathcal{R}_{\mathrm{pop}}(\cdot)$ denotes the population reconstruction operator. Thus, the conventional two-archive cycle can be summarized as convergence preservation by $C^t$, diversity compensation by $A^t$, and offspring generation through archive collaboration. 

However, when this framework is directly applied to LSMOPs, there still remains three limitations. First, it usually assumes that the archive information is useful once it is available, without explicitly evaluating whether the current archive is reliable enough to guide the search. Second, the reproduction operator is generally applied in the full decision space, which is inefficient when only a small subset of variable groups is influential at a given stage. Third, archive updating can preserve existing solutions, but it does not actively repair under-covered front regions or suppress late-stage archive drift. These limitations motivate TRUST-TAEA, which introduces archive trustworthiness to coordinate variable-grouping sparse search, anchor-probing compensation, and archive stabilization.

\section{Proposed TRUST-TAEA}\label{sec3}
This section presents the overall framework of TRUST-TAEA, which consists of four main components: (1) an archive trustworthiness-guided mechanism; (2) an archive trustworthiness-guided variable-grouping sparse search mechanism; (3) an archive trustworthiness-guided anchor-probing compensatory search mechanism; and (4) a checkpoint-based archive stabilization mechanism.

\subsection{Trustworthiness-guided overall framework of TRUST-TAEA}
Building on the decision-variable decomposition described in Section~2.1 and the conventional two-archive framework introduced in Section~2.2, TRUST-TAEA is designed as a trustworthiness-guided two-archive evolutionary framework for LSMOPs. Its main idea is to transform the two archives from passive solution containers into active search regulators. Specifically, the convergence archive $C^t$ is used not only to preserve the convergence-oriented solutions, but also to estimate whether the current archive information is reliable enough to guide the search. The preference archive $A^t$ complements $C^t$ by maintaining distribution information and under-explored regions.

Let the population, convergence archive, and preference archive at generation $t$ be denoted by
\begin{flalign}\label{13}
	& P^t=\left\{\boldsymbol{x}_1^t,\boldsymbol{x}_2^t,\ldots,\boldsymbol{x}_N^t\right\}, \;
	C^t, A^t \subset \Omega, &
\end{flalign}
where $N$ is the population size and $|C^t|=|A^t|=N$. In contrast to conventional TAEAs, TRUST-TAEA first evaluates the reliability of the convergence archive before using it to guide offspring generation. The archive trustworthiness $\mathcal{T}^t$ is estimated from the evolutionary progress and the maturity of the convergence archive, and its detailed definition is given in Section~3.2. A larger $\mathcal{T}^t$ indicates that the archive is more reliable and can provide stronger guidance for exploitation and structural repair; a smaller $\mathcal{T}^t$ implies that the search should maintain more exploration and compensation.

Let $\mathcal{V}$ denote the variable-structure information used in TRUST-TAEA, including the variable groups and the structural targets required for repair. At generation $t$, TRUST-TAEA activates only a subset of informative variable groups from $\mathcal{V}$:
\begin{flalign}\label{14}
	& \mathcal{G}_a^t \subseteq \mathcal{G}, \;
	\mathcal{I}^t=\bigcup_{G_k\in\mathcal{G}_a^t}G_k, &
\end{flalign}
where $\mathcal{G}$ is the variable-group set specified in Section~3.3 and $\mathcal{I}^t$ is the active decision-variable index set.

Accordingly, the offspring-generation process of TRUST-TAEA is formulated as
\begin{flalign}\label{15}
	& \mathcal{Q}^t=
	\mathcal{R}\left(P^t,C^t,A^t,\mathcal{T}^t,\mathcal{G}_a^t\right), &
\end{flalign}
where $\mathcal{R}(\cdot)$ is the trustworthiness-guided reproduction operator. More specifically, the offspring set consists of two complementary parts:
\begin{flalign}\label{16}
	& \mathcal{Q}^t=\mathcal{Q}_{\mathrm{sgs}}^t\cup
	\mathcal{Q}_{\mathrm{probe}}^t, &
\end{flalign}
here, $\mathcal{Q}_{\mathrm{sgs}}^t$ is generated by the variable-grouping sparse search mechanism, which improves search efficiency by modifying only structurally important variables. In contrast, $\mathcal{Q}_{\mathrm{probe}}^t$ is generated by the anchor-probing compensatory search mechanism, which injects candidate solutions into insufficiently covered regions when archive maturity or front coverage is unsatisfactory.

After trustworthiness-guided offspring generation, TRU\-ST-TAEA follows the two-archive update cycle defined in Section~2.2, but with one key difference: the convergence archive is stabilized before it is used to update the preference archive. Specifically, the candidate pool is first constructed as
\begin{flalign}\label{17}
	& U^t=\mathcal{D}\left(P^t\cup \mathcal{Q}^t\cup C^t\cup A^t\right), &
\end{flalign}
where $\mathcal{D}(\cdot)$ removes duplicate solutions. The intermediate convergence archive is obtained by
\begin{flalign}\label{18}
	& \bar{C}^{t+1}=\mathcal{S}_C(U^t,N). &
\end{flalign}

Then, a checkpoint-based stabilization operator is applied:
\begin{flalign}\label{19}
	& C^{t+1}=\mathcal{S}_{\mathrm{ckpt}}\left(\bar{C}^{t+1},C_{\mathrm{ckpt}}^t,p^t\right), &
\end{flalign}
where $C_{\mathrm{ckpt}}^t$ denotes the checkpoint archive and $p^t$ is the normalized evolutionary progress.

Given the stabilized convergence archive, the preference archive and population are updated by the generic operators $\mathcal{S}_A(\cdot)$ and $\mathcal{R}_{\mathrm{pop}}(\cdot)$ defined in Eqs.~(\ref{11}) and (\ref{12}).

Therefore, TRUST-TAEA retains the basic two-archive cycle, while modifying the information flow before archive updating through trustworthiness-guided reproduction and checkpoint stabilization. This control variable coordinates three key behaviors: variable-grouping sparse search, anchor-probing compensation, and checkpoint-based archive stabilization. In this way, the proposed framework explicitly connects archive reliability, variable structure, and search resource allocation, making it more suitable for large-scale multi-objective optimization. The corresponding pseudocode is given in Algorithm~\ref{alg:1}.
\begin{breakablealgorithm}
	\caption{Overall Framework of TRUST-TAEA for LSMOPs}
	\label{alg:1}
	\begin{algorithmic}[1]
		\Require population size $N$, maximum generation $T$, problem $F(\boldsymbol{x})$, variable-structure information $\mathcal{V}$
		\Ensure convergence archive $C^T$, preference archive $A^T$
		
		\State Initialize population $P^0$ and evaluate objective values
		\State Initialize convergence archive $C^0$ and preference archive $A^0$
		\State Initialize checkpoint archive $C_{\mathrm{ckpt}}^0$
		
		\For{$t = 0$ to $T-1$}
		\State Compute archive trustworthiness $\mathcal{T}^t$ by Algorithm~\ref{alg:2}
		\State Generate sparse-search offspring $\mathcal{Q}_{\mathrm{sgs}}^t$ using $\mathcal{T}^t$ and $\mathcal{V}$ by Algorithm~\ref{alg:3}
		\State Generate compensation offspring $\mathcal{Q}_{\mathrm{probe}}^t$ using $\mathcal{T}^t$ and $\mathcal{V}$ by Algorithm~\ref{alg:4}
		\State $\mathcal{Q}^t=\mathcal{Q}_{\mathrm{sgs}}^t\cup\mathcal{Q}_{\mathrm{probe}}^t$
		\State Evaluate $\mathcal{Q}^t$
		\State $U^t=\mathcal{D}\left(P^t\cup\mathcal{Q}^t\cup C^t\cup A^t\right)$
		\State Obtain intermediate convergence archive by Eq.~(\ref{18})
		\State Stabilize convergence archive $C^{t+1}$ by Algorithm~\ref{alg:5}
		\State Update preference archive by Eq.~(\ref{11})
		\State Rebuild population by Eq.~(\ref{12})
		\EndFor
		
		\State \Return $C^T,A^T$
	\end{algorithmic}
\end{breakablealgorithm}

\subsection{Archive trustworthiness guidance mechanism}
The key issue in TRUST-TAEA is how to determine whether the current convergence archive $C^t$ is reliable enough to guide the search. In conventional two-archive algorithms, archive information is usually used directly once it is available. However, in LSMOPs, an archive in early or unstable stages may be sparse, locally collapsed, or structurally biased. Therefore, TRUST-TAEA introduces archive trustworthiness as a control variable.

Let $T$ be the maximum number of generations. The normalized evolutionary progress is
\begin{flalign}\label{20}
	& p^t=\frac{t}{T-1}. &
\end{flalign}

To describe the gradual transition from exploration to exploitation, a stage factor is defined as
\begin{flalign}\label{21}
	& \phi^t=\operatorname{clip}\left(\frac{p^t-\tau_s}{\tau_e-\tau_s},0,1\right), &
\end{flalign}
where $\tau_s$ and $\tau_e$ denote the start and end positions of the transition stage.

The maturity of the convergence archive is evaluated from its nondominated subset $\mathrm{ND}(C^t)$. First, the size maturity is defined as
\begin{flalign}\label{22}
	& M_{\mathrm{size}}^t=\min\left(\frac{|\mathrm{ND}(C^t)|}{\mu N},1\right), &
\end{flalign}
where $N$ is the archive capacity and $\mu$ is a proportional parameter. Second, the coverage maturity is computed according to the occupied bins or reference directions in the normalized objective space:
\begin{flalign}\label{23}
	& M_{\mathrm{cov}}^t=\frac{|\mathcal{B}(C^t)|}{B}, &
\end{flalign}
where $\mathcal{B}(C^t)$ is the set of occupied bins or reference directions and $B$ is the total number of bins or directions. Third, a shape-related maturity term is used to penalize fragmented or collapsed fronts:
\begin{flalign}\label{24}
	& M_{\mathrm{shape}}^t=\frac{1}{1+\kappa(K_{\mathrm{seg}}^t-1)}, &
\end{flalign}
where $K_{\mathrm{seg}}^t$ is the number of detected front segments and $\kappa>0$ is a penalty coefficient. In this study, $\kappa$ is set to 1, so that each additional detected segment directly reduces the shape maturity in a reciprocal form.

The overall convergence-archive maturity is
\begin{flalign}\label{25}
	& \mathcal{M}^t=
	\alpha M_{\mathrm{size}}^t+
	\beta M_{\mathrm{cov}}^t+
	\gamma M_{\mathrm{shape}}^t,\;
	\alpha+\beta+\gamma=1. &
\end{flalign}

Then, archive trustworthiness is defined as
\begin{flalign}\label{26}
	& \mathcal{T}^t=\phi^t\mathcal{M}^t, &
\end{flalign}
\begin{figure}
	\centering
	\includegraphics[width=\linewidth]{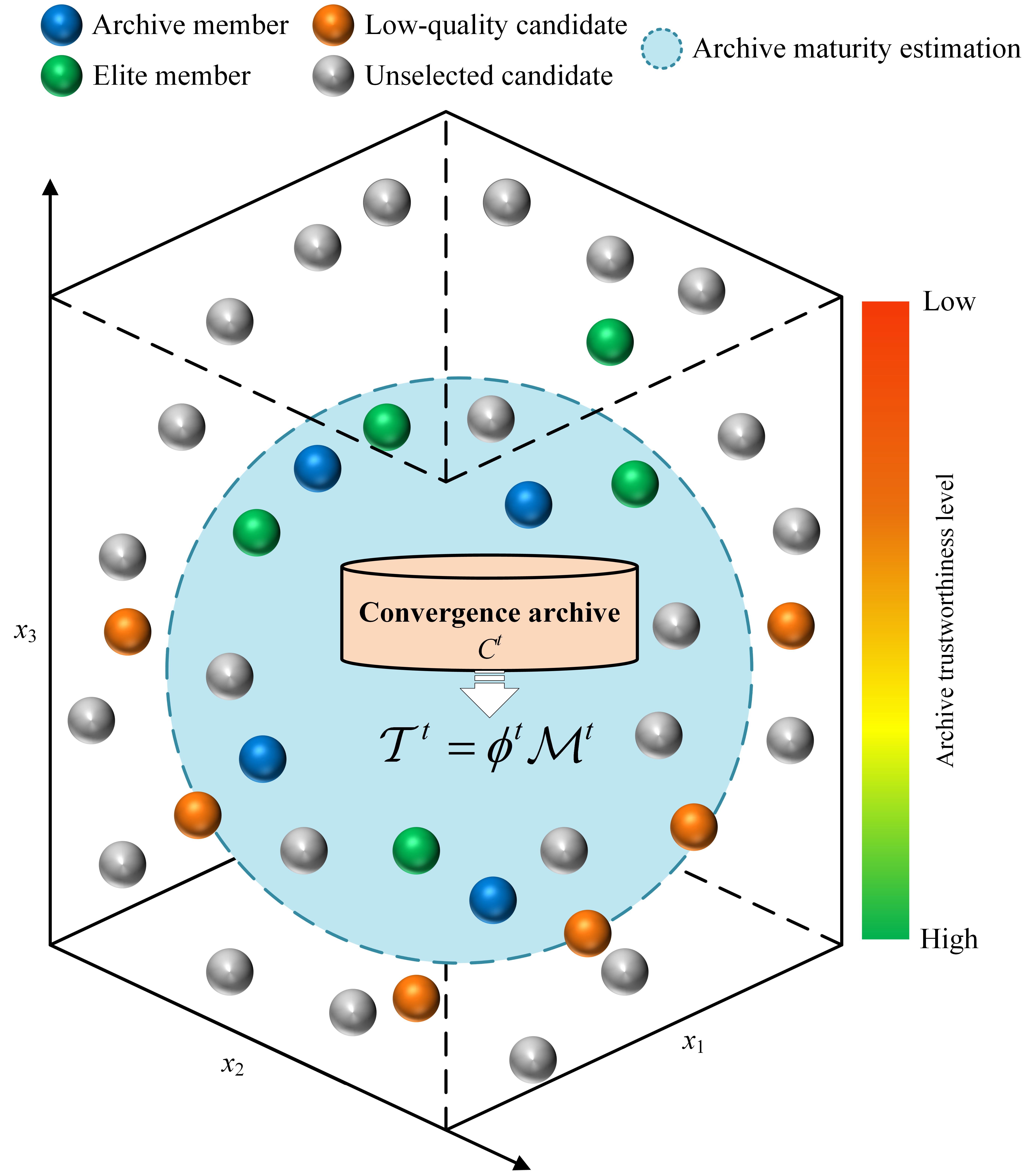} 
	\caption{Illustration of the archive trustworthiness guidance mechanism. The reliability of the convergence archive $C^t$ is estimated by combining the evolutionary stage factor $\phi^t$ and archive maturity $\mathcal{M}^t$, yielding the archive trustworthiness $\mathcal{T}^t=\phi^t\mathcal{M}^t$ to guide subsequent search control.}
	\label{fig2}
\end{figure}
where both $\phi^t\in[0,1]$ and $\mathcal{M}^t\in[0,1]$ are scalar indicators. Thus, $\mathcal{T}^t\in[0,1]$ measures the overall reliability of the convergence archive. This multiplicative form ensures that the archive can be highly trusted only when both the evolutionary stage and archive maturity are sufficiently high. The trustworthiness value is further mapped to the control parameters of offspring generation (Fig.~\ref{fig2}):
\begin{flalign}\label{27}
	& P_{\mathrm{explore}}^t=
	P_{\max}-(P_{\max}-P_{\min})(\mathcal{T}^t)^{\lambda_{\mathrm{exp}}}, &
\end{flalign}
\begin{flalign}\label{28}
	& K_a^t=\left\lceil K_{\min}+(K_{\max}-K_{\min})\mathcal{T}^t\right\rceil, &
\end{flalign}
\begin{flalign}\label{29}
	& \rho^t=\rho_{\min}+(\rho_{\max}-\rho_{\min})\mathcal{T}^t, &
\end{flalign}
here, $P_{\mathrm{explore}}^t$, $K_a^t$, and $\rho^t$ denote the exploration probability, the number of active variable groups, and the structural repair strength, respectively. The parameters satisfy $0\leq P_{\min}<P_{\max}\leq1$, $1\leq K_{\min}\leq K_a^t\leq K_{\max}\leq K$, $0\leq \rho_{\min}<\rho_{\max}\leq1$, and $\lambda_{\mathrm{exp}}>0$. The exponent $\lambda_{\mathrm{exp}}$ controls the nonlinear response of exploration to archive trustworthiness: a larger value keeps exploration higher when $\mathcal{T}^t$ is small, whereas a smaller value makes the algorithm shift earlier toward exploitation. Thus, $\mathcal{T}^t$ provides a unified scheduler that links archive reliability with search behavior. The pseudocode of this module is given in Algorithm~\ref{alg:2}.
\begin{breakablealgorithm}
	\caption{Archive Trustworthiness Guidance Mechanism}
	\label{alg:2}
	\begin{algorithmic}[1]
		\Require generation index $t$, maximum generation $T$, convergence archive $C^t$
		\Ensure archive trustworthiness $\mathcal{T}^t$, exploration probability $P_{\mathrm{explore}}^t$, active group number $K_a^t$, repair strength $\rho^t$
		
		\State Compute normalized progress $p^t$ by Eq. (\ref{20})
		\State Compute stage factor $\phi^t$ by Eq. (\ref{21})
		\State Compute archive size maturity $M_{\mathrm{size}}^t$ by Eq. (\ref{22})
		\State Compute archive coverage maturity $M_{\mathrm{cov}}^t$ by Eq. (\ref{23})
		\State Compute archive shape maturity $M_{\mathrm{shape}}^t$ by Eq. (\ref{24})
		\State Compute archive maturity $\mathcal{M}^t$ by Eq. (\ref{25})
		\State Compute archive trustworthiness $\mathcal{T}^t$ by Eq. (\ref{26})
		\State Compute exploration probability $P_{\mathrm{explore}}^t$ by Eq. (\ref{27})
		\State Compute active group number $K_a^t$ by Eq. (\ref{28})
		\State Compute repair strength $\rho^t$ by Eq. (\ref{29})
		
		\State \Return $\mathcal{T}^t, P_{\mathrm{explore}}^t, K_a^t, \rho^t$
	\end{algorithmic}
\end{breakablealgorithm}

\subsection{Trustworthiness-guided variable-grouping sparse search}
\begin{figure*}[t]
	\centering
	\includegraphics[width=0.8\linewidth]{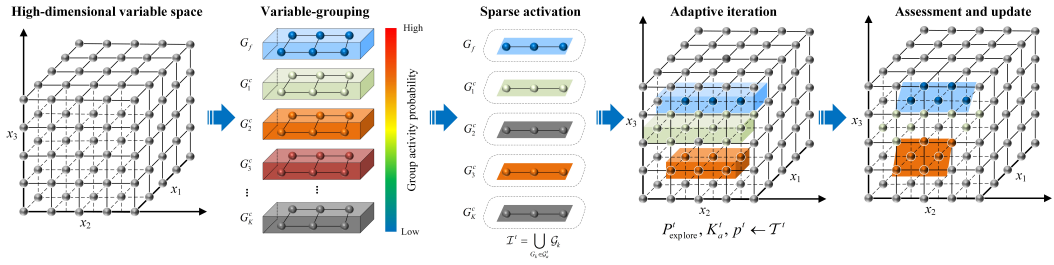}
	\caption{Illustration of the trustworthiness-guided variable-grouping sparse search. Decision variables are decomposed into $G_f$ and $\{G_1^c,\ldots,G_{K_c}^c\}$, and archive trustworthiness $\mathcal{T}^t$ adaptively activates informative groups to form the sparse search subspace $\mathcal{I}^t$ for offspring generation, repair, and update.}
	\label{fig3}
\end{figure*}

The variable-grouping sparse search is the main exploitation branch of TRUST-TAEA. Its purpose is to avoid full-dimensional blind perturbation and concentrate search resources on informative variable groups. This design directly follows the variable-decomposition motivation discussed in Section~2.1. the main idea of this mechanism is in Fig.~\ref{fig3}.

Let the variable groups denoted by
\begin{flalign}\label{30}
	& \mathcal{G}=\{G_f,G_1^c,G_2^c,\ldots,G_{K_c}^c\}, &
\end{flalign}
where $G_f$ contains front-related variables and $G_k^c$ denotes the $k$-th convergence-related group. At generation $t$, an elite subset $E^t\subseteq C^t$ is selected from the convergence archive.

For each group $G_k\in\mathcal{G}$, its search value is estimated by combining dispersion and structural residual. The dispersion term is
\begin{flalign}\label{31}
	& \mathrm{Spr}_k^t=
	\frac{1}{|G_k|}
	\sum_{j\in G_k}
	\frac{\operatorname{std}_{\boldsymbol{x}\in E^t}(x_j)}{u_j-l_j}, &
\end{flalign}
where $u_j$ and $l_j$ are the upper and lower bounds of variable $x_j$. The structural residual is
\begin{flalign}\label{32}
	& \mathrm{Res}_k^t=
	\frac{1}{|E^t||G_k|}
	\sum_{\boldsymbol{x}\in E^t}
	\sum_{j\in G_k}
	|z_j(\boldsymbol{x})-z_k^*|, &
\end{flalign}
where $z_j(\boldsymbol{x})$ is the transformed structural variable and $z_k^*$ is the target value of group $G_k$.

The activity weight of each group is then defined as
\begin{flalign}\label{33}
	& w_k^t=\omega_0+
	\omega_1\widehat{\mathrm{Spr}}_k^t+
	\omega_2\widehat{\mathrm{Res}}_k^t, &
\end{flalign}
where $\omega_0>0$, $\omega_1\geq0$, $\omega_2\geq0$, and $\omega_0+\omega_1+\omega_2=1$. Its sampling probability is
\begin{flalign}\label{34}
	& \pi_k^t=\frac{w_k^t}{\sum_{G_l\in\mathcal{G}}w_l^t}. &
\end{flalign}

According to $\{\pi_k^t\}$, $K_a^t$ groups are selected to form $\mathcal{G}_a^t$, and the active decision-variable set is Eq.~(\ref{15}).

Within the active subspace, a dual-mode reproduction operator is used to generate a mutant vector:
\begin{flalign}\label{35}
	& \boldsymbol{v}_i^t=
	\begin{cases}
		\boldsymbol{a}_i^t+F_i(\boldsymbol{b}_i^t-\boldsymbol{c}_i^t),
		& r_i<P_{\mathrm{explore}}^t,\\
		\boldsymbol{x}_i^t+
		F_i(\boldsymbol{g}_i^t-\boldsymbol{x}_i^t)+
		\lambda F_i(\boldsymbol{b}_i^t-\boldsymbol{c}_i^t),
		& \text{otherwise},
	\end{cases} &
\end{flalign}
where $\boldsymbol{a}_i^t$, $\boldsymbol{b}_i^t$, and $\boldsymbol{c}_i^t$ are selected from $P^t\cup C^t\cup A^t$, $\boldsymbol{g}_i^t$ is an elite-guided vector constructed from $E^t$, $r_i$ is a random number uniformly sampled from $[0,1]$, $F_i$ is the scaling factor, and $\lambda$ is a nonnegative coefficient for the differential perturbation.

Although $\boldsymbol{v}_i^t$ is written in vector form, only the active variables are allowed to inherit its components. The sparse crossover is defined as
\begin{flalign}\label{36}
	& u_{ij}^t=
	\begin{cases}
		v_{ij}^t, & j\in\mathcal{I}^t\cap\mathcal{C}_i^t,\\
		x_{ij}^t, & \text{otherwise},
	\end{cases} &
\end{flalign}
where $\mathcal{C}_i^t$ is the crossover dimension set. Therefore, variables outside $\mathcal{I}^t$ remain unchanged in this sparse-search branch.

For active convergence-related variables, structural repair is further applied:
\begin{flalign}\label{37}
	& u_{ij}^t=(1-\rho^t)u_{ij}^t+\rho^t x_j^*, \;
	j\in \mathcal{I}^t\cap\left(\bigcup_{k=1}^{K_c}G_k^c\right), &
\end{flalign}
where $x_j^*$ denotes the structure-induced target value of the $j$-th decision variable.

This mechanism converts full-dimensional variation into trustworthiness-guided sparse search, improving search efficiency while preserving structural consistency. The pseudocode of this module is given in Algorithm~\ref{alg:3}.
\begin{breakablealgorithm}
	\caption{Trustworthiness-Guided Variable-Grouping Sparse Search}
	\label{alg:3}
	\begin{algorithmic}[1]
		\Require population $P^t$, convergence archive $C^t$, preference archive $A^t$, trustworthiness-controlled parameters $P_{\mathrm{explore}}^t$, $K_a^t$, $\rho^t$, variable-structure information $\mathcal{V}$
		\Ensure sparse-search offspring set $\mathcal{Q}_{\mathrm{sgs}}^t$
		
		\State Obtain variable groups $\mathcal{G}=\{G_f,G_1^c,\ldots,G_{K_c}^c\}$ from $\mathcal{V}$
		\State Select elite subset $E^t \subseteq C^t$
		
		\For{each group $G_k \in \mathcal{G}$}
		\State Compute group dispersion $\mathrm{Spr}_k^t$ by Eq. (\ref{31})
		\State Compute structural residual $\mathrm{Res}_k^t$ by Eq. (\ref{32})
		\State Compute group activity weight $w_k^t$ by Eq. (\ref{33})
		\EndFor
		
		\State Normalize $\{w_k^t\}$ into sampling probabilities $\{\pi_k^t\}$ by Eq. (\ref{34})
		\State Sample $K_a^t$ active groups to form $\mathcal{G}_a^t$
		\State Construct active dimension set $\mathcal{I}^t=\bigcup_{G_k\in\mathcal{G}_a^t}G_k$
		\State $\mathcal{Q}_{\mathrm{sgs}}^t=\emptyset$
		
		\For{each parent solution $\boldsymbol{x}_i^t \in P^t$}
		\If{$\operatorname{rand} < P_{\mathrm{explore}}^t$}
		\State Generate exploration mutant $\boldsymbol{v}_i^t$ by Eq. (\ref{35})
		\Else
		\State Construct elite-guided vector $\boldsymbol{g}_i^t$ from $E^t$
		\State Generate guided mutant $\boldsymbol{v}_i^t$ by Eq. (\ref{35})
		\EndIf
		
		\State Initialize trial solution $\boldsymbol{u}_i^t=\boldsymbol{x}_i^t$
		\For{each dimension $j\in\mathcal{I}^t$}
		\If{$j\in\mathcal{C}_i^t$}
		\State $u_{ij}^t=v_{ij}^t$
		\EndIf
		\EndFor
		
		\For{each active convergence-related dimension $j\in \mathcal{I}^t\cap\left(\bigcup_{k=1}^{K_c}G_k^c\right)$}
		\State Compute structure-induced target $x_j^*$
		\State Repair $u_{ij}^t$ by Eq. (\ref{37})
		\EndFor
		
		\State Add $\boldsymbol{u}_i^t$ into $\mathcal{Q}_{\mathrm{sgs}}^t$
		\EndFor
		
		\State \Return $\mathcal{Q}_{\mathrm{sgs}}^t$
	\end{algorithmic}
\end{breakablealgorithm}

\subsection{Trustworthiness-guided anchor-probing compensatory search}
\begin{figure*}[t]
	\centering
	\includegraphics[width=0.8\linewidth]{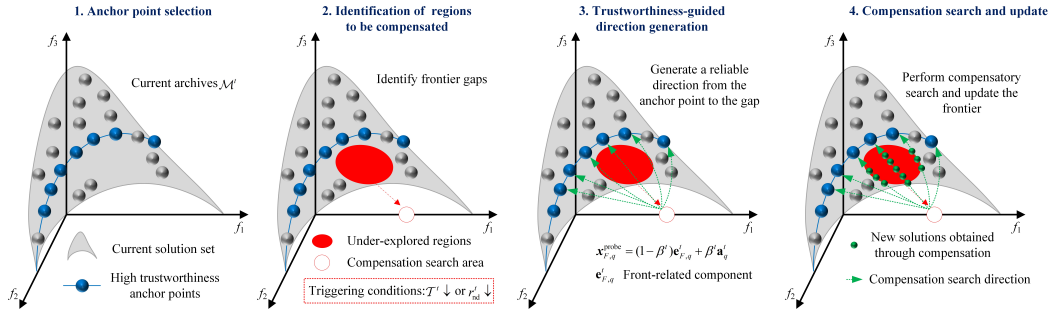}
	\caption{Illustration of the trustworthiness-guided anchor-probing compensatory search. High-trust anchor points are selected from the current archives, and under-explored front regions are identified when archive trustworthiness or the nondominated ratio is low. Probe solutions are then generated along trustworthiness-guided directions to compensate frontier gaps and update the solution set.}
	\label{fig4}
\end{figure*}

Variable-grouping sparse search improves the exploitation efficiency, however, it may still be insufficient when the current front is locally collapsed or poorly covered. Therefore, TRUST-TAEA introduces an anchor-probing compensatory branch to actively supplement under-explored regions. As shown in Fig.~\ref{fig4}.

The compensation intensity is controlled by archive trustworthiness and the nondominated ratio of $C^t$:
\begin{flalign}\label{38}
	& \tilde{\delta}^t =
	\delta_0+\delta_1(1-\mathcal{T}^t)
	+\delta_2(1-r_{\mathrm{nd}}^t), &
\end{flalign}
\begin{flalign}\label{39}
	& \delta^t =
	\begin{cases}
		0, & p^t<p_{\mathrm{start}},\\
		\operatorname{clip}(\tilde{\delta}^t,0,\delta_{\max}),
		& \text{otherwise},
	\end{cases} &
\end{flalign}
where $p_{\mathrm{start}}$ is the activation threshold, $r_{\mathrm{nd}}^t=|\mathrm{ND}(C^t)|/|C^t|$, and $\delta_0,\delta_1,\delta_2\geq0$ are weighting coefficients. The upper bound $\delta_{\max}\in[0,1]$ prevents excessive compensation. The number of probe solutions is $N_{\mathrm{probe}}^t=\lceil N\delta^t\rceil$.

For each probe solution, an elite solution $\boldsymbol{e}_q^t\in E^t$ is selected as the base. Let $\boldsymbol{e}_{F,q}^t$ and $\boldsymbol{e}_{C,q}^t$ denote its front-related and convergence-related components, respectively. Given an anchor point $\boldsymbol{a}_q^t$ in the front-variable space, the front-related component is generated as
\begin{flalign}\label{40}
	& \boldsymbol{x}_{F,q}^{\mathrm{probe}}=
	(1-\beta^t)\boldsymbol{e}_{F,q}^t+
	\beta^t\boldsymbol{a}_q^t, &
\end{flalign}
where $\beta^t$ is the anchor-expansion coefficient.

The convergence-related component is repaired toward the structure-induced target:
\begin{flalign}\label{41}
	& \boldsymbol{x}_{C,q}^{\mathrm{probe}}=
	(1-\rho_{\mathrm{probe}}^t)\boldsymbol{x}_{C,q}^{\mathrm{raw}}+
	\rho_{\mathrm{probe}}^t\boldsymbol{x}_{C,q}^{*}. &
\end{flalign}

Thus, the complete probe solution is
\begin{flalign}\label{42}
	& \boldsymbol{x}_q^{\mathrm{probe}}=
	\left(
	\boldsymbol{x}_{F,q}^{\mathrm{probe}},
	\boldsymbol{x}_{C,q}^{\mathrm{probe}}
	\right), &
\end{flalign}
and the probe offspring set is
\begin{flalign}\label{43}
	& \mathcal{Q}_{\mathrm{probe}}^t=
	\{\boldsymbol{x}_q^{\mathrm{probe}}\}_{q=1}^{N_{\mathrm{probe}}^t}. &
\end{flalign}

This branch complements sparse exploitation. When $C^t$ is not sufficiently trustworthy or its nondominated ratio is low, anchor probing increases the directed exploration and improves the front coverage. The pseudocode of this module is given in Algorithm~\ref{alg:4}. 
\begin{breakablealgorithm}
	\caption{Trustworthiness-Guided Anchor-Probing Compensatory Search}
	\label{alg:4}
	\begin{algorithmic}[1]
		\Require convergence archive $C^t$, archive trustworthiness $\mathcal{T}^t$, evolutionary progress $p^t$, variable-structure information $\mathcal{V}$
		\Ensure compensation offspring set $\mathcal{Q}_{\mathrm{probe}}^t$
		
		\If{$p^t < p_{\mathrm{start}}$}
		\State \Return $\emptyset$
		\EndIf
		
		\State Compute nondominated ratio $r_{\mathrm{nd}}^t=|\mathrm{ND}(C^t)|/|C^t|$
		\State Compute compensation intensity $\delta^t$ by Eq. (\ref{39})
		\State Compute the number of probe solutions $N_{\mathrm{probe}}^t$
		
		\If{$N_{\mathrm{probe}}^t \leq 0$}
		\State \Return $\emptyset$
		\EndIf
		
		\State Construct anchor set $\mathcal{A}^t$ in the front-variable space according to $\mathcal{V}$
		\State $\mathcal{Q}_{\mathrm{probe}}^t=\emptyset$
		\State Select elite subset $E^t\subseteq C^t$
		\For{$q=1$ to $N_{\mathrm{probe}}^t$}
		\State Select elite base solution $\boldsymbol{e}_q^t\in E^t$
		\State Select anchor point $\boldsymbol{a}_q^t\in\mathcal{A}^t$
		\State Generate front-related component $\boldsymbol{x}_{F,q}^{\mathrm{probe}}$ by Eq. (\ref{40})
		\State Generate raw convergence-related component $\boldsymbol{x}_{C,q}^{\mathrm{raw}}$
		\State Construct structure-induced target $\boldsymbol{x}_{C,q}^{*}$ according to $\mathcal{V}$
		\State Repair convergence-related component $\boldsymbol{x}_{C,q}^{\mathrm{probe}}$ by Eq. (\ref{41})
		\State Form probe solution $\boldsymbol{x}_q^{\mathrm{probe}}$ by Eq. (\ref{42})
		\State Add $\boldsymbol{x}_q^{\mathrm{probe}}$ into $\mathcal{Q}_{\mathrm{probe}}^t$
		\EndFor
		
		\State \Return $\mathcal{Q}_{\mathrm{probe}}^t$
	\end{algorithmic}
\end{breakablealgorithm}

\subsection{Checkpoint-based archive stabilization}
In late-stage evolution, the convergence archive may drift away from a previously reliable structural state, especially in high-dimensional problems with strong linkage relationships. To suppress this degradation, TRUST-TAEA introduces a checkpoint-based archive stabilization mechanism.

For each solution $\boldsymbol{x}\in \bar{C}^{t+1}$, the structural residual is defined as
\begin{flalign}\label{44}
	& r(\boldsymbol{x})=
	\frac{1}{D_c}
	\left\|
	\boldsymbol{z}_c(\boldsymbol{x})-\boldsymbol{z}_c^*
	\right\|_1, &
\end{flalign}
where $D_c$ is the number of convergence-related variables, $\boldsymbol{z}_c(\boldsymbol{x})$ is the transformed structural representation, and $\boldsymbol{z}_c^*$ is the structural target. The mean residual of the intermediate convergence archive is
\begin{flalign}\label{45}
	& \bar r^{t+1}=
	\frac{1}{|\bar{C}^{t+1}|}
	\sum_{\boldsymbol{x}\in \bar{C}^{t+1}}r(\boldsymbol{x}). &
\end{flalign}

A composite archive score is used to evaluate the current archive state:
\begin{flalign}\label{46}
	& \Gamma^{t+1}=
	\bar r^{t+1}
	+\lambda_d d^{t+1}
	-\lambda_c \chi^{t+1}
	-\lambda_n r_{\mathrm{nd}}^{t+1}, &
\end{flalign}
where $d^{t+1}$ is the mean normalized objective norm, $\chi^{t+1}$ is the coverage score, and $r_{\mathrm{nd}}^{t+1}$ is the nondominated ratio. A smaller $\Gamma^{t+1}$ indicates a better current archive state.

Let $C_{\mathrm{ckpt}}^t$ denote the saved checkpoint archive, with score 
$\Gamma_{\mathrm{ckpt}}^t$ and mean residual $\bar r_{\mathrm{ckpt}}^t$. 
After the intermediate convergence archive $\bar C^{t+1}$ is obtained, its 
quality is compared with the saved checkpoint. If $\bar C^{t+1}$ shows a 
better archive state, the checkpoint should be refreshed to preserve the 
newly obtained reliable structure. Specifically, the checkpoint is updated when
\begin{flalign}\label{47}
	& \Gamma^{t+1}<\eta_{\Gamma}\Gamma_{\mathrm{ckpt}}^t
	\; \text{or} \;
	\bar r^{t+1}<\eta_r\bar r_{\mathrm{ckpt}}^t, &
\end{flalign}
where $\eta_{\Gamma}<1$ and $\eta_r<1$ are update thresholds.

However, in the later stage of evolution, the current convergence archive may
move away from this reliable checkpoint due to excessive local exploitation or
unstable high-dimensional perturbations. Therefore, after the rollback
mechanism is activated, TRUST-TAEA checks whether the current archive has
significantly deteriorated compared with the saved checkpoint:
\begin{flalign}\label{48}
	& p^t>\tau_b,\;
	\bar r^{t+1}>\gamma_r\bar r_{\mathrm{ckpt}}^t,\;
	\Gamma^{t+1}>\gamma_{\Gamma}\Gamma_{\mathrm{ckpt}}^t, &
\end{flalign}
where $\tau_b$ is the rollback activation threshold, and $\gamma_r>1$ and
$\gamma_{\Gamma}>1$ are degradation thresholds.

If the above condition is satisfied, the current convergence archive is
regarded as degraded and is replaced by the saved checkpoint; otherwise, the
intermediate archive is retained:
\begin{flalign}\label{49}
	& C^{t+1}=
	\begin{cases}
		C_{\mathrm{ckpt}}^t, & \text{if Eq.~(\ref{48}) is satisfied},\\
		\bar{C}^{t+1}, & \text{otherwise}.
	\end{cases} &
\end{flalign}

This mechanism preserves a historically reliable convergence backbone and prevents the two-archive search from losing structural stability in later generations. The pseudocode of this module is given in Algorithm~\ref{alg:5}.
\begin{breakablealgorithm}
	\caption{Checkpoint-Based Archive Stabilization}
	\label{alg:5}
	\begin{algorithmic}[1]
		\Require intermediate convergence archive $\bar{C}^{t+1}$, checkpoint archive $C_{\mathrm{ckpt}}^t$, progress $p^t$
		\Ensure stabilized convergence archive $C^{t+1}$, updated checkpoint archive $C_{\mathrm{ckpt}}^{t+1}$
		
		\State Compute mean structural residual $\bar r^{t+1}$ by Eq. (\ref{45})
		\State Compute composite archive score $\Gamma^{t+1}$ by Eq. (\ref{46})
		\State Compute checkpoint quality $(\bar r_{\mathrm{ckpt}}^t,\Gamma_{\mathrm{ckpt}}^t)$
		
		\If{(\ref{47}) is satisfied}
		\State $C_{\mathrm{ckpt}}^{t+1}=\bar{C}^{t+1}$
		\Else
		\State $C_{\mathrm{ckpt}}^{t+1}=C_{\mathrm{ckpt}}^t$
		\EndIf
		
		\If{(\ref{48}) is satisfied}
		\State $C^{t+1}=C_{\mathrm{ckpt}}^t$
		\Else
		\State $C^{t+1}=\bar{C}^{t+1}$
		\EndIf
		
		\State \Return $C^{t+1}, C_{\mathrm{ckpt}}^{t+1}$
	\end{algorithmic}
\end{breakablealgorithm}

\subsection{Complexity analysis}
The computational complexity of TRUST-TAEA mainly comes from three components: archive-quality estimation, trustworthiness-guided offspring generation, and two-archive environmental selection. Among them, the dominant cost is still the archive update on the union set $U^t$, which is mainly determined by nondominated sorting and archive maintenance. This cost can be expressed as $O(M|U^t|^2)$, where $M$ is the number of objectives.

The additional cost introduced by TRUST-TAEA is mainly associated with the statistical analysis of structural variable groups and the generation of anchor-probe solutions. If the size of the elite subset is denoted by $|E^t|$, the corresponding structural analysis cost is approximately $O(|E^t|D)$, where $D$ is the number of decision variables. For offspring generation, sparse search is performed only on the active variable subset with size $D_a^t$, and thus its cost is $O(ND_a^t)$, with $D_a^t \ll D$ in most generations. The compensatory search introduces an additional cost of approximately $O(N_{\mathrm{probe}}^tD)$, where $N_{\mathrm{probe}}^t$ is the number of probe solutions. Therefore, the per-generation complexity of TRUST-TAEA can be summarized as $ O\left( {M|{U^t}{|^2} + |{E^t}|D + ND_a^t + } \right. $ $\left. {N_{{\rm{probe}}}^tD} \right)$. 

Compared with conventional TAEA, whose main computational cost also comes from two-archive environmental selection and can be roughly expressed as $O(M|U^t|^2)$, TRUST-TAEA does not change the dominant complexity order of the
archive update. The additional terms, including $O(|E^t|D)$ and $O(N_{\mathrm{probe}}^tD)$, are linear with respect to the decision dimension. Moreover, since $|E^t|\leq N$, $N_{\mathrm{probe}}^t$ is controlled by the
compensation intensity, and sparse search is performed on $D_a^t$ active variables with $D_a^t \ll D$ in most generations, the extra computational overhead is controllable. Therefore, TRUST-TAEA preserves the same dominant
complexity as the conventional two-archive framework while introducing only moderate linear costs for trustworthiness-guided search control, making it computationally feasible for LSMOPs. These linear costs are acceptable because they are used to reduce ineffective full-dimensional search and to improve the utilization of limited function evaluations in large-scale decision spaces.

\section{Experiment and results}\label{sec4}
To systematically evaluate the performance of the proposed TRUST-TAEA on LSMOPs, comparative experiments are conducted on the standard LSMOP benchmark suite against several representative algorithms proposed in recent years. The experiments are designed from three perspectives. First, quantitative performance metrics are used to compare the overall performance of different algorithms under different numbers of objectives and decision-space dimensions. Second, Pareto front visualizations are used to examine the differences among algorithms in terms of front approximation and distribution maintenance. Third, parameter sensitivity analysis is performed to investigate the influence of key parameter settings on the performance of TRUST-TAEA. Through these experiments, the effectiveness of the proposed algorithm is evaluated from multiple aspects, including convergence, distribution, and robustness. All numerical experiments were conducted on a workstation equipped with a 12th Gen Intel Core i5-12490F processor at 3.00 GHz and 16 GB RAM.

\subsection{Benchmark problems and performance metrics}
The LSMOP benchmark suite effectively captures the challenges caused by high-dimensional decision spaces, complex variable interactions, and diverse Pareto front shapes in large-scale multi-objective optimization. It is therefore widely used to evaluate the optimization capability of algorithms in large-scale scenarios. To comprehensively examine algorithmic performance under different problem scales, experiments are conducted on LSMOP benchmark instances with two and three objectives, where the number of decision variables is set to 500, 1000, 2000, and 5000.
\begin{table*}[htbp]
	\caption{HV results and statistical comparisons of TRUST-TAEA and the competing algorithms on 2-, 3-objective LSMOPs with different decision dimensions.}\label{Tab1}
	\centering
	\resizebox{\textwidth}{!}{%
		\begin{tabular}{lcccccccccc}
			\hline
			\multirow{2.5}{*}{Problem} & \multirow{2.5}{*}{M} & \multirow{2.5}{*}{D} & \multicolumn{8}{@{}c@{}}{Algorithm}  \\ \cmidrule{4-11}%
			&  &  &  DGEA & FLEA & PCPSO & LERD & NNCSO & AMRCSO  & MSCSO & TRUST-TAEA \\
			\midrule
			\multirow{8}{*}{LSMOP1} & \multirow{4}{*}{2} & 500 & 1.4345E-1(9.44E-2)- & 1.1376E-1(2.85E-2)- & 1.2346E-2(1.21E-2)- & 1.5391E-3(3.21E-3)- & 1.4976E-1(6.32E-2)- & 2.5591E-1(3.56E-3)- & 5.9663E-1(3.94E-3)- &\textbf{ 6.9088E-1(1.44E-3)}  \\ 
			& & 1000 & 1.3645E-1(7.33E-2)- & 1.3125E-1(2.82E-2)- & 0.0000E+0(0.00E+0)- & 0.0000E+0(0.00E+0)- & 2.0733E-1(6.76E-2)- & 2.4908E-1(1.11E-3)- & 3.6374E-1(5.92E-3)- & \textbf{6.9301E-1(9.20E-4)}  \\ 
			& & 2000 & 1.4498E-1(8.24E-2)- & 1.1528E-1(1.13E-2)- & 0.0000E+0(0.00E+0)- & 0.0000E+0(0.00E+0)- & 1.1947E-1(5.94E-2)- & 2.7959E-1(1.19E-3)- & 3.2123E-1(3.91E-3)- & \textbf{6.9316E-1(8.72E-4)} \\
			& & 5000 & 2.1241E-1(7.23E-2)- & 1.1465E-1(3.05E-2)- & 0.0000E+0(0.00E+0)- & 0.0000E+0(0.00E+0)- & 9.7707E-2(6.05E-2)- & 2.8341E-1(5.24E-4)- & 4.5427E-1(5.12E-3)- & \textbf{6.9415E-1(1.52E-3)} \\
			& \multirow{4}{*}{3} & 500 & 1.6332E-1(9.57E-2)- & 1.2332E-1(4.52E-2)- & 2.3877E-4(3.94E-4)- & 9.1947E-4(1.56E-3)- & 1.4423E-1(3.34E-2)- & 4.4553E-1(2.55E-3)- & 6.5327E-1(2.88E-2)- & \textbf{7.2652E-1(2.57E-2)} \\
			& & 1000 & 1.0914E-1(1.10E-1)- & 1.1014E-1(3.14E-2)- & 0.0000E+0(0.00E+0)- & 9.4710E-5(3.22E-4)- & 1.3325E-1(4.00E-2)- & 3.9438E-1(1.23E-3)- & 4.0277E-1(3.98E-2)- & \textbf{7.4220E-1(6.14E-2)} \\
			& & 2000 & 6.7128E-2(8.35E-2)- & 5.9086E-2(5.27E-2)- & 0.0000E+0(0.00E+0)- & 0.0000E+0(0.00E+0)- & 1.2450E-1(3.03E-2)- & 3.5919E-1(2.21E-3)- & 4.1780E-1(5.49E-2)- & \textbf{7.3258E-1(3.16E-2)} \\
			& & 5000 & 1.2486E-1(8.27E-2)- & 2.8813E-2(6.27E-2)- & 0.0000E+0(0.00E+0)- & 9.1122E-6(4.72E-5)- & 1.0351E-1(3.65E-2)- & 3.9651E-1(2.58E-3)- & 4.0254E-1(5.15E-2)- & \textbf{7.3115E-1(4.21E-2)} \\
			\multirow{8}{*}{LSMOP2} & \multirow{4}{*}{2} & 500 & 4.6789E-1(1.15E-3)- & 6.5270E-1(3.62E-3)- & 4.6125E-1(7.40E-4)- & 6.1058E-1(1.31E-2)- & 5.4536E-1(7.82E-4)- & 5.5316E-1(2.94E-4)- & 4.6382E-1(1.39E-4)- & \textbf{6.9509E-1(1.59E-3)}  \\ 
			& & 1000 & 5.9475E-1(8.36E-4)- & 6.0653E-1(8.45E-3)- & 4.7294E-1(3.46E-4)- & 4.5458E-1(3.99E-3)- & 6.3905E-1(5.48E-4)- & 6.4979E-1(8.60E-5)- & 4.6758E-1(2.27E-4)- & \textbf{6.9605E-1(1.27E-3)}  \\
			& & 2000 & 5.4731E-1(3.41E-4)- & 5.4185E-1(2.08E-3)- & 4.8313E-1(1.97E-4)- & 5.9398E-1(8.67E-4)- & 5.9593E-1(2.47E-4)- & 5.7310E-1(3.15E-5)- & 6.8252E-1(5.66E-4)- & \textbf{6.9634E-1(1.21E-3)}  \\ 
			& & 5000 & 5.9488E-1(1.07E-3)- & 6.5074E-1(3.57E-4)- & 5.0465E-1(9.26E-5)- & 5.0825E-1(3.70E-4)- & 5.5303E-1(3.20E-4)- & 6.1314E-1(2.83E-5)- & \textbf{6.9504E-1(1.03E-4)}$ \approx $ & \textbf{6.9645E-01(1.71E-3)}  \\ 
			& \multirow{4}{*}{3} & 500 & 8.9305E-1(4.30E-3)+ & 6.6358E-1(3.34E-3)- & 6.5552E-1(2.17E-3)- & 6.9901E-1(2.38E-3)- & 8.5040E-1(9.25E-4)+ & 8.8000E-1(5.14E-4)+ & \textbf{9.7424E-1(3.84E-4)+} & 7.2389E-1(9.41E-3)  \\ 
			& & 1000 & 9.3436E-1(1.53E-3)+ & 7.7089E-1(2.69E-3)+ & 8.7898E-1(1.01E-3)+ & 9.3516E-1(4.52E-4)+ & 7.2186E-1(5.16E-4)$ \approx $ & 7.4892E-1(2.56E-4)+ & \textbf{9.7053E-1(7.06E-4)+} & 7.2778E-1(6.96E-3)  \\ 
			& & 2000 & 6.7269E-1(4.36E-4)- & 7.5915E-1(4.21E-3)+ & \textbf{9.4013E-1(2.03E-3)+} & 9.2499E-1(2.47E-4)+ & 6.8279E-1(4.26E-4)- & 7.6891E-1(4.59E-4)+ & 7.1594E-1(2.25E-4)- & 7.3439E-1(9.10E-3)  \\ 
			& & 5000 & 6.8101E-1(3.83E-4)- & 9.6845E-1(6.14E-3)+ & \textbf{9.8000E-1(7.92E-4)+} & 7.1412E-1(3.66E-4)- & 8.3515E-1(2.81E-4)+ & 8.2729E-1(2.80E-4)+ & 9.2481E-1(3.82E-4)+ & 7.3536E-1(7.04E-3)  \\ 
			\multirow{8}{*}{LSMOP3} & \multirow{4}{*}{2} & 500 & 0.0000E+0(0.00E+0)- & 0.0000E+0(0.00E+0)- & 0.0000E+0(0.00E+0)- & 0.0000E+0(0.00E+0)- & 0.0000E+0(0.00E+0)- & 0.0000E+0(0.00E+0)- & 0.0000E+0(0.00E+0)- & \textbf{1.0715E-2(1.99E-2)}  \\ 
			& & 1000 & 0.0000E+0(0.00E+0)- & 0.0000E+0(0.00E+0)- & 0.0000E+0(0.00E+0)- & 0.0000E+0(0.00E+0)- & 0.0000E+0(0.00E+0)- & 0.0000E+0(0.00E+0)- & 0.0000E+0(0.00E+0)- & \textbf{ 1.8432E-2(2.82E-2)}  \\ 
			& & 2000 & 0.0000E+0(0.00E+0)- & 0.0000E+0(0.00E+0)- & 0.0000E+0(0.00E+0)- & 0.0000E+0(0.00E+0)- & 0.0000E+0(0.00E+0)- & 0.0000E+0(0.00E+0)- & 0.0000E+0(0.00E+0)- & \textbf{2.2162E-2(3.81E-2)}  \\ 
			& & 5000 & 0.0000E+0(0.00E+0)- & 0.0000E+0(0.00E+0)- & 0.0000E+0(0.00E+0)- & 0.0000E+0(0.00E+0)- & 0.0000E+0(0.00E+0)- & 0.0000E+0(0.00E+0)- & 0.0000E+0(0.00E+0)- & \textbf{2.8494E-2(4.74E-2)}  \\
			& \multirow{4}{*}{3} & 500 & 6.5053E-2(4.01E-2)- & 8.2530E-2(2.43E-2)- & 0.0000E+0(0.00E+0)- & 0.0000E+0(0.00E+0)- & 7.0319E-2(3.45E-2)- & 1.1970E-2(2.77E-2)- & 2.5894E-2(3.53E-2)- & \textbf{9.0612E-1(7.53E-2)}  \\ 
			& & 1000 & 7.9527E-2(2.78E-2)- & 2.4008E-2(4.19E-2)- & 0.0000E+0(0.00E+0)- & 0.0000E+0(0.00E+0)- & 7.0194E-2(3.63E-2)- & 4.7555E-2(5.18E-2)- & 3.1500E-2(3.02E-2)- & \textbf{8.9067E-1(8.28E-2)}  \\ 
			& & 2000 & 7.0529E-2(3.71E-2)- & 0.0000E+0(0.00E+0)- & 0.0000E+0(0.00E+0)- & 0.0000E+0(0.00E+0)- & 4.5221E-2(5.69E-2)- & 3.4878E-2(5.10E-2)- & 1.2491E-2(3.93E-2)- & \textbf{9.0251E-1(5.74E-2)}  \\
			& & 5000 & 4.0821E-2(5.55E-2)- & 0.0000E+0(0.00E+0)- & 0.0000E+0(0.00E+0)- & 0.0000E+0(0.00E+0)- & 4.2773E-2(3.56E-2)- & 1.7057E-2(2.54E-2)- & 4.0415E-2(5.94E-2)- & \textbf{9.2956E-1(6.26E-2)}  \\ 
			\multirow{8}{*}{LSMOP4} & \multirow{4}{*}{2} & 500 & 6.1653E-1(5.10E-3)- & 6.0791E-1(4.65E-3)- & 5.0857E-1(2.40E-3)- & 5.1995E-1(6.18E-3)- & 5.1051E-1(1.83E-3)- & 4.6052E-1(1.25E-3)- & 4.8811E-1(1.17E-4)- & \textbf{6.8569E-1(1.36E-3)}  \\ 
			& & 1000 & 6.4884E-1(2.30E-3)- & 5.1809E-1(5.73E-3)- & 5.2687E-1(1.24E-3)- & 4.6766E-1(6.01E-3)- & 5.0720E-1(1.34E-3)- & 5.6769E-1(3.02E-4)- & 4.2836E-1(3.63E-4)- & \textbf{6.9101E-1(1.21E-3)}  \\ 
			& & 2000 & 6.2291E-1(1.25E-3)- & 5.4095E-1(3.34E-3)- & 5.2221E-1(2.46E-4)- & 5.2812E-1(2.77E-3)- & 6.5636E-1(7.75E-4)- & 4.5005E-1(7.97E-5)- & 4.8270E-1(7.93E-5)- & \textbf{6.9285E-1(1.23E-3)}  \\ 
			& & 5000 & 6.0267E-1(4.71E-4)- & 4.8086E-1(1.35E-3)- & 5.2231E-1(4.93E-5)- & 6.2192E-1(5.15E-4)- & 4.6882E-1(2.40E-4)- & 5.2812E-1(3.39E-5)- & 4.4280E-1(5.34E-4)- & \textbf{6.9560E-1(1.32E-3)}  \\ 
			& \multirow{4}{*}{3} & 500 & \textbf{7.6661E-1(1.17E-2)+} & 7.1420E-1(1.28E-2)- & 7.1586E-1(2.85E-3)- & 6.0192E-1(7.35E-3)- & 6.4204E-1(4.24E-3)- & 6.3973E-1(1.73E-3)- & 7.2571E-1(2.38E-4)- & 7.4550E-1(2.61E-2)  \\ 
			& & 1000 & 7.2971E-1(3.39E-3)+ & \textbf{8.6589E-1(5.37E-3)+} & 8.1261E-1(2.88E-3)+ & 7.2480E-1(3.32E-3)$ \approx $ & 6.7161E-1(1.51E-3)- & 8.3378E-1(5.74E-4)+ & 6.1355E-1(2.11E-3)- & 7.1710E-1(1.33E-2)  \\ 
			& & 2000 & 9.5174E-1(2.90E-3)+ & 7.9558E-1(2.92E-3)+ & 7.3590E-1(1.55E-3)+ & 8.5378E-1(1.07E-3)+ & 9.5663E-1(1.24E-3)+ & 7.9686E-1(9.17E-4)+ & \textbf{9.6742E-1(8.29E-4)+ }& 7.2466E-1(1.22E-2)  \\ 
			& & 5000 & \textbf{9.5773E-1(1.73E-3)+} & 9.3677E-1(3.14E-3)+ & 8.4223E-1(1.60E-3)+ & 7.1262E-1(5.63E-4)- & 7.3474E-1(2.61E-4)+ & 7.1230E-1(4.55E-4)- & 6.6513E-1(3.72E-4)- & 7.2706E-1(9.14E-3)  \\ 
			\multirow{8}{*}{LSMOP5} & \multirow{4}{*}{2} & 500 & 9.1355E-3(2.60E-2)- & 9.2932E-2(4.21E-17)- & 0.0000E+0(0.00E+0)- & 5.4087E-3(2.65E-2)- & 3.2247E-2(3.42E-2)- & 1.0005E-1(3.73E-17)- & 1.2173E-2(3.01E-2)- & \textbf{3.8360E-1(7.16E-3)}  \\ 
			& & 1000 & 0.0000E+0(0.00E+0)- & 8.5695E-2(3.55E-17)- & 0.0000E+0(0.00E+0)- & 0.0000E+0(0.00E+0)- & 2.0053E-2(2.84E-2)- & 8.3371E-2(5.24E-17)- & 1.1522E-1(1.45E-2)- & \textbf{3.8774E-1(5.40E-3)}  \\ 
			& & 2000 & 7.8746E-3(2.32E-2)- & 6.7582E-2(3.63E-2)- & 0.0000E+0(0.00E+0)- & 0.0000E+0(0.00E+0)- & 1.0541E-3(2.41E-3)- & 7.5107E-2(4.99E-17)- & 6.7737E-2(5.05E-2)- & \textbf{3.8739E-1(5.96E-3)}  \\ 
			& & 5000 & 0.0000E+0(0.00E+0)- & 5.7492E-2(5.01E-2)- & 0.0000E+0(0.00E+0)- & 0.0000E+0(0.00E+0)- & 8.0656E-3(2.21E-2)- & 8.6512E-2(4.21E-17)- & 5.1913E-2(4.68E-2)- & \textbf{3.9239E-1(4.80E-3)}  \\ 
			& \multirow{4}{*}{3} & 500 & 2.6978E-1(1.11E-1)- & 1.03919E-1(4.89E-2)- & 0.0000E+0(0.00E+0)- & 1.6201E-2(7.64E-2)- & 3.3662E-1(1.01E-1)- & 3.8514E-1(2.22E-3)- & 2.5638E-1(5.05E-2)- & \textbf{5.4264E-1(1.03E-1)}  \\ 
			& & 1000 & 1.8468E-1(1.33E-1)- & 1.0936E-1(5.21E-2)- & 0.0000E+0(0.00E+0)- & 0.0000E+0(0.00E+0)- & 2.5152E-1(1.43E-1)- & 3.2969E-1(4.75E-3)- & 3.4185E-1(8.39E-2)- & \textbf{5.0046E-1(1.47E-1)}  \\ 
			& & 2000 & 1.5541E-1(1.48E-1)- & 7.3948E-2(7.44E-2)- & 0.0000E+0(0.00E+0)- & 0.0000E+0(0.00E+0)- & 1.9124E-1(1.50E-1)- & 3.2001E-1(6.54E-3)- & 2.0111E-1(9.16E-2)- & \textbf{5.0720E-1(1.25E-1)}  \\ 
			& & 5000 & 1.8075E-1(1.58E-1)- & 2.1324E-2(2.94E-2)- & 0.0000E+0(0.00E+0)- & 0.0000E+0(0.00E+0)- & 2.1452E-1(1.47E-1)- & 3.2110E-1(6.21E-3)- & 2.3431E-1(1.32E-1)- & \textbf{5.3483E-1(1.46E-1)}  \\
			\multirow{8}{*}{LSMOP6} & \multirow{4}{*}{2} & 500 & 5.0234E-2(2.45E-2)- & 4.5679E-2(2.78E-2)- & 0.0000E+0(0.00E+0)- & 6.7371E-3(7.13E-3)- & 4.7579E-3(8.80E-3)- & 4.2628E-2(2.85E-2)- & 9.2044E-2(6.12E-6)- & \textbf{1.7023E-1(1.09E-2)}  \\ 
			& & 1000 & 6.5949E-2(1.72E-2)- & 8.0667E-2(2.23E-2)- & 2.6002E-2(1.88E-2)- & 4.1477E-2(2.97E-3)- & 3.0582E-2(8.22E-3)- & 7.5018E-2(1.23E-3)- & 4.2765E-2(6.18E-3)- & \textbf{1.7743E-1(1.37E-2)}  \\ 
			& & 2000 & 6.7425E-2(1.29E-2)- & 6.5508E-2(7.61E-3)- & 3.5263E-2(3.35E-2)- & 6.0212E-2(1.98E-3)- & 5.6495E-2(2.13E-2)- & 8.0105E-2(6.58E-3)- & 9.7321E-2(4.86E-7)- & \textbf{1.7922E-1(1.32E-2)}  \\ 
			& & 5000 & 9.7354E-2(7.63E-3)- & 7.8861E-2(2.98E-3)- & 5.7622E-2(4.76E-2)- & 8.7440E-2(7.50E-4)- & 9.0521E-2(9.29E-4)- & 9.7105E-2(3.12E-4)- & 9.5622E-2(3.67E-4)- & \textbf{1.9066E-1(1.92E-2)}  \\ 
			& \multirow{4}{*}{3} & 500 & 9.2328E-3(1.61E-2)- & 5.5549E-4(2.50E-3)- & 0.0000E+0(0.00E+0)- & 0.0000E+0(0.00E+0)- & 0.0000E+0(0.00E+0)- & 1.5941E-2(4.98E-3)- & 5.9178E-3(6.91E-3)- & \textbf{9.4585E-1(9.56E-2)}  \\ 
			& & 1000 & 1.8381E-2(1.57E-2)- & 0.0000E+0(0.00E+0)- & 0.0000E+0(0.00E+0)- & 0.0000E+0(0.00E+0)- & 0.0000E+0(0.00E+0)- & 2.2633E-2(2.16E-3)- & 1.2286E-2(1.27E-2)- & \textbf{9.5447E-1(7.75E-2)}  \\ 
			& & 2000 & 1.0142E-2(8.53E-3)- & 0.0000E+0(0.00E+0)- & 0.0000E+0(0.00E+0)- & 0.0000E+0(0.00E+0)- & 0.0000E+0(0.00E+0)- & 1.7056E-2(2.11E-3)- & 4.8865E-3(5.87E-3)- & \textbf{9.6656E-1(5.83E-2)}  \\ 
			& & 5000 & 7.3301E-3(7.79E-3)- & 0.0000E+0(0.00E+0)- & 0.0000E+0(0.00E+0)- & 0.0000E+0(0.00E+0)- & 0.0000E+0(0.00E+0)- & 2.0864E-2(4.12E-3)- & 9.4145E-3(1.33E-2)- & \textbf{8.8666E-1(1.23E-1)}  \\ 
			\multirow{8}{*}{LSMOP7} & \multirow{4}{*}{2} & 500 & 0.0000E+0(0.00E+0) & 0.0000E+0(0.00E+0) & 0.0000E+0(0.00E+0) & 0.0000E+0(0.00E+0) & 0.0000E+0(0.00E+0) & 0.0000E+0(0.00E+0) & 0.0000E+0(0.00E+0) & 0.0000E+0(0.00E+0)  \\ 
			& & 1000 & 0.0000E+0(0.00E+0)- & 0.0000E+0(0.00E+0)- & 0.0000E+0(0.00E+0)- & 0.0000E+0(0.00E+0)- & 0.0000E+0(0.00E+0)- & 0.0000E+0(0.00E+0)- & 0.0000E+0(0.00E+0)- & \textbf{1.9920E-3(8.68E-3)}  \\ 
			& & 2000 & 0.0000E+0(0.00E+0) & 0.0000E+0(0.00E+0) & 0.0000E+0(0.00E+0) & 0.0000E+0(0.00E+0) & 0.0000E+0(0.00E+0) & 0.0000E+0(0.00E+0) & 0.0000E+0(0.00E+0) & 0.0000E+0(0.00E+0)  \\ 
			& & 5000 & 0.0000E+0(0.00E+0)- & 0.0000E+0(0.00E+0)- & 0.0000E+0(0.00E+0)- & 0.0000E+0(0.00E+0)- & 0.0000E+0(0.00E+0)- & 0.0000E+0(0.00E+0)- & 0.0000E+0(0.00E+0)- & \textbf{1.8024E-3(3.84E-3)}  \\ 
			& \multirow{4}{*}{3} & 500 & 0.0000E+0(0.00E+0)- & 0.0000E+0(0.00E+0)- & 0.0000E+0(0.00E+0)- & 5.3572E-4(2.34E-3)- & 0.0000E+0(0.00E+0)- & 0.0000E+0(0.00E+0)- & 6.7237E-2(4.93E-2)- & \textbf{5.9996E-1(8.36E-2)}  \\ 
			& & 1000 & 6.7324E-3(5.95E-3)- & 0.0000E+0(0.00E+0) & 0.0000E+0(0.00E+0) & 4.7754E-4(2.58E-3)- & 9.7642E-4(2.40E-3)- & 0.0000E+0(0.00E+0)- & 3.6382E-2(1.53E-3)- & \textbf{5.6063E-1(7.62E-2)}  \\ 
			& & 2000 & 4.9214E-2(6.58E-3)- & 1.2062E-2(8.12E-3)- & 0.0000E+0(0.00E+0)- & 1.6835E-2(9.75E-3)- & 4.2945E-2(2.37E-3)- & 1.9943E-2(0.00E+0)- & 8.9164E-2(1.16E-3)- & \textbf{5.4120E-1(1.09E-1)}  \\ 
			& & 5000 & 6.8721E-2(1.23E-3)- & 4.7639E-2(2.01E-3)- & 1.4003E-2(2.06E-2)- & 7.4951E-2(2.65E-3)- & 6.7476E-2(1.11E-3)- & 7.2338E-2(3.62E-17)- & 6.4004E-2(2.90E-3)- & \textbf{5.7421E-1(5.75E-2)}  \\ 
			\multirow{8}{*}{LSMOP8} & \multirow{4}{*}{2} & 500 & 3.6069E-3(1.89E-2)- & 7.4393E-2(2.39E-3)- & 0.0000E+0(0.00E+0)- & 1.6186E-2(3.76E-2)- & 8.9735E-2(1.806E-2)- & 1.2013E-1(6.80E-3)- & 8.3006E-2(2.73E-3)- & \textbf{4.0081E-1(3.00E-3)}  \\ 
			& & 1000 & 0.0000E+0(0.00E+0)- & 8.3195E-2(5.23E-17)- & 0.0000E+0(0.00E+0)- & 0.0000E+0(0.00E+0)- & 5.0490E-2(4.45E-2)- & 9.3462E-2(3.91E-17)- & 2.1590E-1(6.77E-3)- & \textbf{4.0388E-1(2.22E-3)}  \\ 
			& & 2000 & 1.0006E-2(3.54E-2)- & 8.3488E-2(2.48E-2)- & 0.0000E+0(0.00E+0)- & 0.0000E+0(0.00E+0)- & 1.9634E-2(4.83E-2)- & 8.9771E-2(4.89E-17)- & 1.0088E-1(1.94E-2)- & \textbf{4.0501E-1(3.45E-3)}  \\ 
			& & 5000 & 4.8769E-3(2.51E-2)- & 8.2974E-2(3.66E-2)- & 0.0000E+0(0.00E+0)- & 0.0000E+0(0.00E+0)- & 1.2711E-2(3.27E-2)- & 9.3802E-2(4.37E-17)- & 1.0541E-1(7.94E-3)- & \textbf{4.1310E-1(2.47E-3)}  \\
			& \multirow{4}{*}{3} & 500 & 4.6084E-1(1.70E-2)+ & 7.3158E-2(6.25E-2)- & 1.7003E-2(2.25E-2)- & 3.9076E-1(1.98E-2)+ & 2.3356E-1(1.71E-1)- & \textbf{5.3922E-1(3.60E-3)+} & 3.0723E-1(8.67E-2)- & 3.6271E-1(2.10E-2)  \\ 
			& & 1000 & 2.9587E-1(1.64E-1)- & 8.2356E-2(9.59E-4)- & 4.3443E-2(3.24E-2)- & 2.9607E-1(1.34E-1)- & 2.3779E-1(1.79E-1)- & \textbf{4.4854E-1(3.50E-3)+} & 3.8769E-1(2.44E-2)+ & 3.4883E-1(8.08E-3)  \\ 
			& & 2000 & 3.2545E-1(1.06E-1)- & 7.3871E-2(1.23E-3)- & 3.8595E-2(4.17E-2)- & 1.3020E-1(7.42E-2)- & 1.6642E-1(1.60E-1)- & \textbf{4.3546E-1(1.35E-3)+} & 2.8824E-1(1.16E-1)- & 3.4912E-1(8.14E-3)  \\ 
			& & 5000 & 2.0484E-1(2.19E-1)- & 7.6926E-2(1.12E-3)- & 3.9423E-2(5.13E-2)- & 9.8386E-2(1.22E-2)- & 2.0835E-1(1.13E-1)- & \textbf{5.5330E-1(1.53E-3)+} & 4.2032E-1(9.40E-2)+ & 3.3993E-1(8.79E-3)  \\ 
			\multirow{8}{*}{LSMOP9} & \multirow{4}{*}{2} & 500 & 0.0000E+0(0.00E+0)- & 9.5598E-2(1.62E-2)- & 4.3019E-2(4.02E-2)- & 4.5311E-2(2.44E-2)- & 8.9931E-2(1.56E-2)- & 1.9277E-1(2.04E-2)- & 9.7052E-2(1.52E-2)- & \textbf{7.2065E-1(7.69E-3)}  \\ 
			& & 1000 & 0.0000E+0(0.00E+0)- & 7.7733E-2(3.05E-2)- & 0.0000E+0(0.00E+0)- & 2.5021E-2(1.86E-2)- & 9.1362E-2(7.95E-3)- & 2.0245E-1(1.99E-3)- & 8.9689E-2(3.46E-2)- & \textbf{7.4837E-1(3.78E-3)}  \\ 
			& & 2000 & 0.0000E+0(0.00E+0)- & 2.2913E-2(4.73E-2)- & 0.0000E+0(0.00E+0)- & 0.0000E+0(0.00E+0)- & 8.6167E-2(6.29E-3)- & 2.0141E-1(1.81E-2)- & 7.6463E-2(3.83E-2)- & \textbf{7.6266E-1(3.27E-3)}  \\ 
			& & 5000 & 0.0000E+0(0.00E+0)- & 0.0000E+0(0.00E+0)- & 2.3817E-5(8.33E-5)- & 0.0000E+0(0.00E+0)- & 1.0144E-1(7.33E-3)- & 1.7712E-1(1.54E-2)- & 7.9679E-2(2.42E-2)- & \textbf{7.8120E-1(2.76E-3)}  \\ 
			& \multirow{4}{*}{3} & 500 & 0.0000E+0(0.00E+0)- & 1.0683E-1(4.49E-2)- & 0.0000E+0(0.00E+0)- & 1.9610E-2(1.32E-2)- & 1.5965E-1(1.51E-2)- & 1.9699E-1(4.25E-3)- & 1.1092E-1(2.11E-2)- & \textbf{2.6979E-1(7.81E-3)}  \\ 
			& & 1000 & 0.0000E+0(0.00E+0)- & 3.6253E-2(6.16E-2)- & 0.0000E+0(0.00E+0)- & 0.0000E+0(0.00E+0)- & 1.6838E-1(7.36E-2)- & 1.9239E-1(6.03E-3)- & 1.1355E-1(1.51E-4)- & \textbf{2.7918E-1(6.48E-3)}  \\
			& & 2000 & 4.8759E-5(7.90E-3)- & 1.3109E-2(4.94E-2)- & 0.0000E+0(0.00E+0)- & 0.0000E+0(0.00E+0)- & 1.7298E-1(4.86E-2)- & 1.7841E-1(3.34E-3)- & 1.5805E-1(3.48E-2)- & \textbf{2.7937E-1(5.59E-3)}  \\ 
			& & 5000 & 0.0000E+0(0.00E+0)- & 1.1659E-3(4.88E-3)- & 0.0000E+0(0.00E+0)- & 0.0000E+0(0.00E+0)- & 1.7118E-1(5.81E-2)- & 1.7578E-1(3.49E-3)- & 1.0897E-1(3.23E-2)- & \textbf{2.8371E-1(8.38E-3)}  \\ \hline
			\multicolumn{3}{@{}c@{}}{+/-/$ \approx $} & 7/65/0 & 6/66/0 & 6/66/0 & 4/67/1 & 4/67/1 & 6/65/1 & 5/66/1 & \\
			\bottomrule 
	\end{tabular}}
\end{table*}

For performance evaluation, the inverted generational distance plus indicator (IGD$^+$) \cite{bib36} and the hypervolume indicator (HV) \cite{bib37} are used to assess the performance of all algorithms. In addition, the Wilcoxon rank-sum test at a significance level of 0.05 is adopted for statistical analysis. The symbols ``$+$'', ``$-$'', and ``$\approx$'' indicate that the compared algorithm performs significantly better than, significantly worse than, or comparably to the proposed algorithm, respectively.
\begin{table*}[htbp]
	\caption{IGD$ ^+ $ results and statistical comparisons of TRUST-TAEA and the competing algorithms on 2-, 3-objective LSMOPs with different decision dimensions.}\label{Tab2}
	\centering
	\resizebox{\textwidth}{!}{%
		\begin{tabular}{lcccccccccc}
			\hline
			\multirow{2.5}{*}{Problem} & \multirow{2.5}{*}{M} & \multirow{2.5}{*}{D} & \multicolumn{8}{@{}c@{}}{Algorithm}  \\ \cmidrule{4-11}%
			&  &  &  DGEA & FLEA & PCPSO & LERD & NNCSO & AMRCSO  & MSCSO & TRUST-TAEA \\
			\midrule
			\multirow{8}{*}{LSMOP1} & \multirow{4}{*}{2} & 500 & 4.5773E-1(1.50E-1)- & 1.0982E+0(7.19E-2)- & 7.9488E-1(5.43E-2)- & 1.0751E+0(1.51E-1)- & 3.6280E-1(9.51E-2)- & 2.3305E-1(5.44E-3)- & 1.7060E-1(4.76E-3)- & \textbf{9.8527E-3(6.39E-4)}  \\ 
			& & 1000 & 4.4146E-1(1.28E-1)- & 6.8818E-1(4.83E-2)- & 1.1560E+0(4.63E-2)- & 1.3443E+0(2.18E-1)- & 3.6973E-1(9.20E-2)- & 2.6515E-1(3.17E-3)- & 1.1609E-1(4.64E-3)- & \textbf{8.8782E-3(4.41E-4)}  \\ 
			& & 2000 & 4.9205E-1(2.42E-1)- & 7.3807E-1(5.22E-2)- & 1.4330E+0(4.35E-2)- & 1.8024E+0(2.47E-1)- & 4.1239E-1(9.72E-2)- & 2.8799E-1(1.68E-3)- & 1.8479E-1(4.31E-3)- & \textbf{8.7943E-3(3.93E-4)} \\
			& & 5000 & 4.7656E-1(1.29E-1)- & 7.8273E-1(1.35E-1)- & 1.8482E+0(7.93E-2)- & 1.5198E+0(1.86E-1)- & 5.2049E-1(1.02E-1)- & 3.6126E-1(7.28E-4)- & 1.8102E-1(5.57E-3)- & \textbf{8.0781E-3(2.94E-4)} \\
			& \multirow{4}{*}{3} & 500 & 6.1140E-1(1.55E-1)- & 7.4682E-1(1.23E-1)- & 1.5084E+0(3.42E-1)- & 1.1754E+0(1.52E-1)- & 7.3894E-1(4.32E-2)- & 3.5777E-1(3.66E-3)- & 2.7709E-1(2.44E-2)- & \textbf{6.3375E-2(3.32E-3)} \\
			& & 1000 & 6.7617E-1(1.42E-1)- & 8.2396E-1(1.46E-1)- & 1.3146E+0(1.34E-1)- & 1.4659E+0(1.11E-1)- & 7.5644E-1(8.09E-2)- & 3.9508E-1(2.26E-3)- & 2.5238E-1(2.03E-2)- & \textbf{5.9125E-2(7.33E-3)} \\
			& & 2000 & 7.8204E-1(1.64E-1)- & 9.8319E-1(7.69E-1)- & 1.6785E+0(8.62E-2)- & 1.5474E+0(1.29E-1)- & 7.5253E-1(6.99E-2)- & 3.7643E-1(2.52E-3)- & 2.9220E-1(3.62E-2)- & \textbf{6.1101E-2(3.26E-3)} \\
			& & 5000 & 6.0571E-1(1.12E-1)- & 1.4665E+0(8.30E-1)- & 1.4217E+0(1.40E-1)- & 1.5940E+0(1.58E-1)- & 7.1561E-1(4.27E-2)- & 3.1956E-1(2.49E-3)- & 2.9050E-1(5.72E-2)- & \textbf{6.2682E-2(3.60E-3)} \\
			\multirow{8}{*}{LSMOP2} & \multirow{4}{*}{2} & 500 & 1.1753E-2(8.59E-4)- & 2.5788E-2(3.64E-3)- & 4.0764E-2(6.29E-4)- & 2.6290E-2(7.74E-3)- & 1.6491E-2(7.43E-4)- & 1.0483E-2(1.60E-4)- & \textbf{4.2392E-3(5.68E-5)+} & 6.8042E-3(6.97E-4)  \\ 
			& & 1000 & 7.2016E-3(5.11E-4)- & 1.7229E-2(5.48E-3)- & 2.0007E-2(2.53E-4)- & 1.8889E-2(3.07E-3)- & 1.0030E-2(3.11E-4)- & 7.0901E-3(4.74E-5)- & \textbf{5.6294E-3(1.29E-4)+} & 6.1568E-3(6.08E-4)  \\
			& & 2000 & 4.8983E-3(3.60E-4)+ & 1.6357E-2(1.78E-3)- & 1.0426E-2(1.66E-4)- & 1.2974E-2(6.39E-4)- & 5.8686E-3(1.01E-4)+ & 4.8283E-3(2.41E-5)+ & \textbf{3.7332E-3(1.39E-5)+} & 6.1925E-3(6.85E-4)  \\ 
			& & 5000 & 4.7897E-3(9.83E-4)+ & 9.8446E-3(2.46E-4)- & 6.7269E-3(6.07E-5)- & 7.0178E-3(3.14E-4)- & \textbf{3.5744E-3(1.88E-4)+} & 3.6330E-3(2.06E-5)+ & 3.9362E-3(3.06E-5)+ & 5.5862E-3(5.80E-4)  \\ 
			& \multirow{4}{*}{3} & 500 & 5.7134E-2(2.00E-3)- & 7.6565E-2(3.36E-3)- & 7.2636E-2(1.38E-3)- & 6.6161E-2(1.52E-3)- & 4.6972E-2(5.71E-4)+ & \textbf{3.6358E-2(6.02E-4)+} & 4.2690E-2(2.08E-4)+ & 4.8321E-2(1.99E-3)  \\ 
			& & 1000 & \textbf{3.5531E-2(1.34E-3)+} & 5.2146E-2(2.56E-3)- & 5.4465E-2(7.92E-4)- & 4.1224E-2(3.70E-4)+ & 4.0601E-2(5.10E-4)+ & 4.1723E-2(3.88E-4)+ & 4.6231E-2(2.99E-4)- & 4.2958E-2(2.00E-3)  \\ 
			& & 2000 & \textbf{3.3316E-2(2.06E-4)+} & 4.8293E-2(3.98E-3)- & 4.7691E-2(1.13E-3)- & 3.7599E-2(1.49E-4)+ & 3.5247E-2(1.98E-4)+ & 3.4980E-2(4.31E-4)+ & 4.1689E-2(8.77E-5)- & 4.0298E-2(1.96E-3)  \\ 
			& & 5000 & 3.2069E-2(1.79E-4)+ & 4.3180E-2(2.39E-3)- & 4.8406E-2(9.86E-4)- & 3.4117E-2(1.29E-4)+ & \textbf{2.8905E-2(1.10E-4)+} & 4.0101E-2(2.65E-4)- & 4.2139E-2(1.24E-4)- & 3.8716E-2(1.65E-3)  \\ 
			\multirow{8}{*}{LSMOP3} & \multirow{4}{*}{2} & 500 & 1.5044E+0(3.27E-1)- & 1.5050E+0(8.81E-4)- & 2.7387E+1(7.80E-1)- & 8.0823E+0(7.58E+0)- & 4.4286E+0(3.17E+0)- & 1.8184E+0(1.90E-3)- & 1.8559E+0(1.23E+0)- & \textbf{2.4100E-1(1.15E-1)}  \\ 
			& & 1000 & 1.8640E+0(2.83E-1)- & 3.0849E+0(4.84E+0)- & 2.6141E+1(7.20E-1)- & 1.4881E+1(2.72E+0)- & 3.2520E+0(2.70E+0)- & 1.5880E+0(1.23E-3)- & 1.4784E+0(3.72E-2)- & \textbf{1.8529E-1(7.27E-2)}  \\ 
			& & 2000 & 3.6012E+0(2.50E+0)- & 9.3006E+0(7.40E+0)- & 2.8340E+1(7.54E-1)- & 1.7768E+1(1.56E+0)- & 2.2005E+0(1.07E+0)- & 1.6298E+0(5.40E-4)- & 1.5903E+0(2.10E-3)- & \textbf{1.9778E-1(5.85E-2)}  \\ 
			& & 5000 & 4.4251E+0(3.41E+0)- & 1.6788E+1(2.92E+0)- & 2.8716E+1(5.43E-1)- & 1.7910E+1(1.12E+0)- & 3.6595E+0(3.34E+0)- & 1.3118E+0(1.30E-4)- & 1.5903E+0(2.52E-3)- & \textbf{1.7313E-1(5.96E-2)}  \\
			& \multirow{4}{*}{3} & 500 & 1.9416E+0(1.98E+0)- & 1.0724E+0(4.42E-1)- & 1.8848E+1(3.60E+0)- & 8.2392E+0(3.50E-1)- & 1.0285E+0(8.28E-1)- & 1.1202E+0(1.22E-1)- & 1.2120E+0(3.13E+0)- & \textbf{2.3314E-1(5.43E-2)}  \\ 
			& & 1000 & 1.0633E+0(9.98E-1)- & 5.1733E+0(3.23E+0)- & 2.1429E+1(3.63E+0)- & 1.1542E+1(3.15E-1)- & 1.1066E+0(6.79E-1)- & 8.5896E-1(4.33E-2)- & 2.7366E+0(1.86E+0)- & \textbf{2.4309E-1(6.74E-2)}  \\ 
			& & 2000 & 1.0926E+0(4.96E-1)- & 1.2267E+1(5.19E+0)- & 2.0627E+1(7.38E+0)- & 9.4861E+0(2.77E-1)- & 1.2474E+0(2.04E+0)- & 1.0931E+0(1.07E-1)- & 5.4879E+0(2.78E+0)- & \textbf{2.4514E-1(6.72E-2)}  \\
			& & 5000 & 2.5035E+0(3.25E+0)- & 1.5068E+1(4.53E+0)- & 3.0557E+1(5.56E+1)- & 1.0764E+1(3.58E-1)- & 1.1080E+0(5.73E-1)  & 1.0595E+0(5.83E-2)- & 3.0792E+0(2.89E+0)- & \textbf{2.1010E-1(5.34E-2)}  \\ 
			\multirow{8}{*}{LSMOP4} & \multirow{4}{*}{2} & 500 & 5.3998E-2(4.93E-3)- & 5.6529E-2(4.36E-3)- & 7.5181E-2(2.10E-3)- & 4.7527E-2(7.41E-3)- & 4.0106E-2(1.50E-3)- & 4.2311E-2(9.14E-4)- & 1.6523E-2(9.76E-5)- & \textbf{1.0811E-2(5.93E-4)}  \\ 
			& & 1000 & 2.6132E-2(1.89E-3)- & 3.5644E-2(3.12E-3)- & 5.3273E-2(6.53E-4)- & 3.3434E-2(6.07E-3)- & 2.9473E-2(1.10E-3)- & 1.8538E-2(1.72E-4)- & 1.9368E-2(3.74E-4)- & \textbf{8.6898E-3(3.58E-4)}  \\ 
			& & 2000 & 1.5875E-2(8.77E-4)- & 2.3450E-2(2.60E-3)- & 2.9226E-2(1.66E-4)- & 1.8126E-2(1.32E-3)- & 1.1614E-2(3.68E-4)- & 9.2917E-3(4.74E-5)- & 4.3356E-3(3.86E-4)+ & \textbf{7.8285E-3(6.23E-4)}  \\ 
			& & 5000 & 6.3040E-3(3.47E-4)+ & 1.2523E-2(1.06E-3)- & 1.5012E-2(6.25E-5)- & 8.5608E-3(5.31E-4)- & \textbf{6.0759E-3(1.61E-4)+} & 5.1288E-3(2.01E-5)+ & 6.2253E-3(6.90E-5)+ & 6.6042E-3(5.40E-4)  \\ 
			& \multirow{4}{*}{3} & 500 & 1.2876E-1(7.01E-3)- & 1.6563E-1(8.78E-3)- & 1.8443E-1(1.61E-3)- & 1.3382E-1(8.46E-3)- & 9.5308E-2(2.47E-3)- & 9.9477E-2(1.66E-3)- & \textbf{4.7698E-2(4.14E-4)+} & 7.4810E-2(3.89E-3)  \\ 
			& & 1000 & 6.9819E-2(2.97E-3)- & 1.2460E-1(4.18E-3)- & 9.9003E-2(1.81E-3)- & 1.1113E-1(3.71E-3)- & \textbf{6.0908E-2(1.16E-3)+} & 5.3472E-2(4.82E-4)+ & 6.2938E-2(1.46E-3)$ \approx $ & 6.2930E-2(2.85E-3)  \\ 
			& & 2000 & 5.5421E-2(2.23E-3)- & 6.9116E-2(2.44E-3)- & 6.5530E-2(1.41E-3)- & 5.3173E-2(1.30E-3)- & 4.8930E-2(7.80E-4)+ & 5.3954E-2(6.67E-4)- & \textbf{4.2886E-2(1.43E-4)+} & 5.2123E-2(2.29E-3)  \\ 
			& & 5000 & 4.1419E-2(1.19E-3)+ & 5.5244E-2(1.92E-3)- & 4.0472E-2(7.56E-4)+ & 3.8100E-2(4.73E-4)+ & 3.6496E-2(1.93E-4)+ & \textbf{3.4412E-2(4.56E-4)+} & 4.4837E-2(2.46E-4)$ \approx $ & 4.4457E-2(1.99E-3)  \\ 
			\multirow{8}{*}{LSMOP5} & \multirow{4}{*}{2} & 500 & 4.9794E+0(1.45E+0)- & 7.9313E-1(1.90E-16)- & 1.8713E+0(1.14E-1)- & 3.0129E+0(5.92E-1)- & 1.4994E+0(9.59E-1)- & 6.4015E-1(1.99E-16)- & 1.7670E-1(7.52E-2)- & \textbf{2.1268E-2(3.58E-3)}  \\ 
			& & 1000 & 6.2526E+0(1.07E+0)- & 1.3608E-1(2.61E-16)- & 2.6091E+0(1.27E-1)- & 4.3120E+0(6.10E-1)- & 1.8114E+0(1.27E+0)- & 7.4536E-1(3.05E-16)- & 2.4057E-1(1.88E-2)- & \textbf{1.9530E-2(2.77E-3)}  \\ 
			& & 2000 & 5.1236E+0(2.29E+0)- & 9.4194E-1(1.96E-1)- & 2.8694E+0(9.33E-2)- & 4.8072E+0(1.08E+0)- & 2.0436E+0(1.04E+0)- & 8.5057E-1(1.88E-16)- & 6.1720E-1(7.85E-2)- & \textbf{1.9648E-2(3.15E-3)}  \\ 
			& & 5000 & 4.9470E+0(1.50E+0)- & 8.6287E-1(2.83E-1)- & 3.1202E+0(8.06E-2)- & 5.0550E+0(8.22E-1)- & 2.4984E+0(8.21E-1)- & 6.3224E-1(2.24E-16)- & 6.2263E-1(6.12E-2)- & \textbf{1.7045E-2(2.43E-3)}  \\ 
			& \multirow{4}{*}{3} & 500 & 6.1622E-1(2.19E-1)- & 1.0043E+0(2.47E-1)- & 2.7092E+0(4.63E-1)- & 2.3979E+0(8.86E-1)- & 5.3815E-1(2.15E-1)- & 5.3375E-1(7.31E-3)- & 4.1012E-1(8.96E-2)- & \textbf{5.8031E-2(5.97E-3)}  \\ 
			& & 1000 & 7.0148E-1(3.89E-1)- & 9.4087E-1(2.65E-1)- & 3.7449E+0(3.23E-1)- & 3.3066E+0(7.86E-1)- & 7.6105E-1(4.32E-1)- & 4.5111E-1(6.69E-3)- & 4.2355E-1(1.04E-1)- & \textbf{5.7390E-2(1.22E-2)}  \\ 
			& & 2000 & 9.3688E-1(4.44E-1)- & 1.0953E+0(4.52E-1)- & 4.1580E+0(3.82E-1)- & 2.9213E+0(3.10E-1)- & 8.1579E-1(2.95E-1)- & 5.4814E-1(8.44E-3)- & 5.2006E-1(5.64E-2)- & \textbf{5.9038E-2(7.38E-3)}  \\ 
			& & 5000 & 9.3837E-1(5.18E-1)- & 1.8637E+0(8.62E-1)- & 4.0481E+0(5.99E-1)- & 3.1374E+0(4.88E-1)- & 6.5792E-1(2.41E-1)- & 4.7242E-1(1.28E-2)- & 5.4370E-1(4.48E-2)- & \textbf{5.8845E-2(9.53E-3)}  \\
			\multirow{8}{*}{LSMOP6} & \multirow{4}{*}{2} & 500 & 6.5739E-1(2.04E-1)- & 4.2535E-1(1.26E-1)- & 6.5117E+2(4.48E+2)- & 8.9717E-1(5.53E-2)- & 8.0526E-1(1.31E-1)- & 3.1551E-1(4.26E-2)- & 6.7350E-1(7.18E-5)- & \textbf{1.2684E-1(2.59E-2)}  \\ 
			& & 1000 & 4.6093E-1(2.44E-1)- & 3.9206E-1(1.14E-1)- & 3.3879E+2(5.15E+2)- & 7.3313E-1(1.08E-1)- & 7.6553E-1(5.57E-1)- & 3.0646E-1(6.25E-2)- & 6.7448E-1(2.79E-2)- & \textbf{1.1837E-1(2.83E-2)}  \\ 
			& & 2000 & 6.1221E-1(1.97E-1)- & 4.4882E-1(1.41E-1)- & 6.7064E+2(2.01E+3)- & 6.4912E-1(1.01E-3)- & 9.7226E-1(7.64E-1)- & 2.6909E-1(1.17E-2)- & 6.7463E-1(1.47E-5)- & \textbf{1.1189E-1(2.50E-2)}  \\ 
			& & 5000 & 6.9272E-1(1.79E-1)- & 6.0856E-1(1.72E-1)- & 8.2519E+2(2.27E+3)- & 8.3199E-1(3.86E-4)- & 7.1863E-1(1.20E-1)- & 2.8732E-1(3.20E-4)- & 6.7538E-1(1.31E-2)- & \textbf{1.0005E-1(2.15E-2)}  \\ 
			& \multirow{4}{*}{3} & 500 & 2.9587E+1(6.28E+1)- & 4.2539E+1(2.33E+2)- & 2.9036E+2(2.67E+2)- & 7.5190E+0(9.20E+0)- & 1.3939E+0(1.85E-1)- & 6.4700E-1(1.15E-2)- & 8.1511E-1(4.02E-1)- & \textbf{6.1705E-2(7.35E-2)}  \\ 
			& & 1000 & 2.9509E+1(6.45E+1)- & 1.5053E+0(2.25E-1)- & 8.5011E+2(8.49E+2)- & 3.4258E+1(2.84E+1)- & 1.5555E+0(1.36E-1)- & 8.4550E-1(2.67E-2)- & 8.4084E-1(3.34E-1)- & \textbf{6.8017E-2(6.47E-2)}  \\ 
			& & 2000 & 8.9620E+1(1.12E+2)- & 9.5881E+1(1.84E+2)- & 1.0690E+3(1.29E+3)- & 2.0117E+2(1.26E+2)- & 1.7831E+0(1.70E-1)- & 6.4467E-1(1.13E-2)- & 1.2243E+0(4.34E-1)- & \textbf{6.0925E-2(5.77E-2)}  \\ 
			& & 5000 & 1.2334E+2(2.77E+2)- & 7.2176E+1(1.87E+2)- & 1.6896E+3(1.45E+3)- & 4.7371E+2(2.08E+2)- & 1.3388E+0(1.28E-1)- & 5.6664E-1(2.98E-2)- & 1.0750E+0(4.26E-1)- & \textbf{1.1144E-1(8.73E-2)}  \\ 
			\multirow{8}{*}{LSMOP7} & \multirow{4}{*}{2} & 500 & 3.3105E+3(1.33E+3)- & 1.3383E+0(2.66E-4)- & 2.8479E+2(5.64E+1)- & 4.0853E+0(5.80E-1)- & 1.7265E+0(3.15E-4)- & 1.3414E+0(9.34E-4)- & 1.5226E+0(1.78E-3)- & \textbf{3.2890E-1(9.38E-2)}  \\ 
			& & 1000 & 2.1635E+3(1.86E+3)- & 4.7051E+0(1.44E+1)- & 9.1706E+2(1.13E+2)- & 7.7490E+0(1.45E+0)- & 1.3771E+0(1.10E-3)- & 1.6673E+0(6.58E-4)- & 2.1257E+0(1.93E+0)- & \textbf{2.7308E-1(1.05E-1)}  \\ 
			& & 2000 & 4.3800E+3(2.66E+3)- & 4.5988E+1(1.42E+2)- & 1.4936E+3(1.69E+2)- & 1.1841E+2(5.47E+1)- & 1.6031E+0(5.03E-3)- & 1.5975E+0(3.99E-4)- & 3.0281E+1(2.67E+2)- & \textbf{3.2173E-1(1.05E-1)}  \\ 
			& & 5000 & 4.2163E+3(3.77E+3)- & 2.5196E+0(3.96E+0)- & 2.1046E+3(1.53E+2)- & 2.4974E+3(7.26E+2)- & 3.9783E+1(1.78E+2)- & 1.4551E+0(2.49E-4)- & 3.2166E+1(1.19E+2)- & \textbf{2.1334E-1(9.90E-2)}  \\ 
			& \multirow{4}{*}{3} & 500 & 1.0481E+0(7.17E-2)- & 8.9101E-1(1.40E-1)- & 1.1616E+4(1.16E+4)- & 9.1211E-1(9.52E-2)- & 1.1497E+0(2.39E-2)- & 8.0293E-1(7.72E-3)- & 8.8546E-1(4.11E-2)- & \textbf{1.5470E-1(2.21E-2)}  \\ 
			& & 1000 & 1.1108E+0(7.75E-2)- & 9.3222E-1(1.45E-2)- & 2.2436E+4(1.00E+4)- & 1.2051E+0(3.27E-2)- & 1.0451E+0(1.51E-2)- & 7.5642E-1(6.80E-3)- & 8.7422E-1(4.20E-2)- & \textbf{1.6258E-1(2.04E-2)}  \\ 
			& & 2000 & 9.2450E-1(6.08E-2)- & 9.1698E-1(8.68E-3)- & 1.6817E+4(1.28E+4)- & 9.7911E-1(6.23E-3)- & 9.4010E-1(4.63E-3)- & 8.1110E-1(3.18E-3)- & 9.1709E-1(2.58E-2)- & \textbf{1.4095E-1(2.97E-2)}  \\ 
			& & 5000 & 8.5375E-1(8.82E-2)- & 8.5998E-1(1.58E-3)- & 1.6723E+4(1.42E+4)- & 8.3056E-1(1.12E-3)- & 1.0548E+0(6.96E-3)- & 7.5408E-1(2.15E-3)- & 8.9483E-1(3.73E-2)- & \textbf{1.0500E-1(2.47E-2)}  \\ 
			\multirow{8}{*}{LSMOP8} & \multirow{4}{*}{2} & 500 & 1.5468E+0(1.05E+0)- & 7.3096E-1(1.07E-3)- & 1.4350E+0(8.36E-2)- & 6.9449E-1(2.00E-1)- & 7.6330E-1(5.51E-3)- & 3.2284E-1(1.93E-2)- & 3.8195E-1(1.54E-2)- & \textbf{1.2052E-2(1.42E-3)}  \\ 
			& & 1000 & 2.5622E+0(1.38E+0)- & 8.5389E-1(2.98E-16)- & 2.0155E+0(7.26E-2)- & 1.5072E+0(3.66E-1)- & 8.7273E-1(3.35E-1)- & 5.4382E-1(1.19E-2)- & 1.8819E-1(9.66E-3)- & \textbf{1.0578E-2(1.23E-3)}  \\ 
			& & 2000 & 2.5019E+0(1.33E+0)- & 8.3285E-1(5.84E-4) & 2.2581E+0(9.32E-2)- & 1.9699E+0(1.08E+0)- & 9.4430E-1(2.20E-1)- & 6.5095E-1(1.60E-16)- & 4.3193E-1(3.54E-2)- & \textbf{9.6154E-3(1.39E-3)}  \\ 
			& & 5000 & 2.2379E+0(1.09E+0)- & 6.7188E-1(9.56E-2)- & 2.8210E+0(8.10E-2)- & 2.4787E+0(5.91E-1)- & 9.9930E-1(2.16E-1)- & 8.6909E-1(2.26E-16)- & 4.3206E-1(2.63E-2)- & \textbf{5.2149E-3(9.94E-4)}  \\
			& \multirow{4}{*}{3} & 500 & 1.6291E-1(2.73E-2)- & 7.4695E-1(1.06E-1)- & 2.1149E+0(1.68E+0)- & 2.4747E-1(4.82E-2)- & 5.1516E-1(4.10E-1)- & 1.0084E-1(1.43E-3)- & 1.6385E-1(4.62E-2)- & \textbf{5.3936E-2(2.45E-3)}  \\ 
			& & 1000 & 3.4939E-1(2.66E-1)- & 6.6539E-1(7.70E-2)- & 2.2075E+0(2.07E+0)- & 3.0005E-1(1.13E-1)- & 5.4885E-1(2.32E-1)- & 9.0809E-2(1.86E-3)- & 1.1783E-1(3.23E-2)- & \textbf{5.2281E-2(1.83E-3)}  \\ 
			& & 2000 & 2.8536E-1(2.50E-1)- & 9.4876E-1(8.81E-2)- & 2.1017E+0(1.93E+0)- & 4.5676E-1(8.88E-2)- & 4.2969E-1(3.15E-1)- & 6.5830E-2(1.09E-3)- & 1.5529E-1(6.46E-2)- & \textbf{5.1008E-2(2.28E-3)}  \\ 
			& & 5000 & 4.0774E-1(3.32E-1)- & 8.4366E-1(8.59E-2)- & 2.7230E+0(2.59E+0)- & 5.4137E-1(4.95E-2)- & 4.3219E-1(2.85E-1)- & 8.3580E-2(1.32E-3)- & 1.8149E-1(8.33E-2)- & \textbf{5.3447E-2(1.42E-3)}  \\ 
			\multirow{8}{*}{LSMOP9} & \multirow{4}{*}{2} & 500 & 5.0219E+0(2.15E+0)- & 7.1565E-1(1.03E-1)- & 8.8165E-1(4.24E-1)- & 8.5961E-1(8.70E-2)- & 8.1599E-1(9.09E-2)- & 1.4551E-1(1.28E-1)- & 7.9026E-1(8.33E-2)- & \textbf{2.5695E-2(2.08E-3)}  \\ 
			& & 1000 & 9.3476E+0(3.24E+0)- & 9.2114E-1(3.94E-1)- & 1.9862E+0(3.74E-1)- & 1.0419E+0(3.82E-1)- & 7.5814E-1(5.11E-2)- & 5.2602E-2(7.07E-3)- & 7.6018E-1(1.31E-1)- & \textbf{1.5573E-2(1.09E-3)}  \\ 
			& & 2000 & 1.2891E+1(5.75E+0)- & 1.8790E+0(9.77E-1)- & 2.1239E+0(4.27E-1)- & 4.7492E+0(1.37E+0)- & 7.5621E-1(3.87E-2)- & 8.1725E-2(7.58E-2)- & 8.1173E-1(2.87E-1)- & \textbf{1.0744E-2(8.19E-4)}  \\ 
			& & 5000 & 1.8055E+1(2.83E+0)- & 3.4398E+0(6.83E-1)- & 1.5257E+0(1.55E-1)- & 2.6880E+0(2.36E+0)- & 8.0318E-1(2.53E-2)- & 5.9978E-2(1.05E-1)- & 7.9635E-1(9.27E-2)- & \textbf{4.8316E-3(5.18E-4)}  \\ 
			& \multirow{4}{*}{3} & 500 & 2.1609E+1(7.09E+0)- & 1.0470E+0(6.02E-1)- & 1.0016E+1(9.69E-1)- & 2.4455E+0(2.47E-1)- & 6.1395E-1(1.44E-1)- & 6.4222E-1(9.27E-3)- & 1.1947E+0(2.66E-1)- & \textbf{9.6319E-2(7.10E-3)}  \\ 
			& & 1000 & 3.5026E+1(1.37E+1)- & 3.6618E+0(1.94E+0)- & 1.4665E+1(6.53E-1)- & 5.3745E+0(4.91E+0)- & 9.2523E-1(1.34E+0)- & 6.3847E-1(8.26E-3)- & 1.1459E+0(7.80E-4)- & \textbf{8.6799E-2(7.32E-3)}  \\
			& & 2000 & 4.6120E+1(1.35E+1)- & 1.2782E+1(7.26E+0)- & 2.3114E+1(8.08E-1)- & 9.8628E+0(3.88E+0)- & 9.7958E-1(2.76E+0)- & 5.0795E-1(5.82E-3)- & 1.1321E+0(2.24E+0)- & \textbf{8.5516E-2(7.64E-3)}  \\ 
			& & 5000 & 4.6422E+1(1.50E+1)- & 2.3231E+1(5.09E+0)- & 2.4097E+1(6.93E-1)- & 1.9570E+1(9.34E+0)- & 9.1378E-1(7.04E-1)- & 6.2433E-1(3.19E-3)- & 1.2886E+0(7.54E-1)- & \textbf{8.1369E-2(8.73E-3)}  \\ \hline
			\multicolumn{3}{@{}c@{}}{+/-/$ \approx $} & 7/65/0 & 0/72/0 & 1/71/0 & 4/68/0 & 10/62/0 & 8/64/0 & 9/61/2 & \\
			\bottomrule 
	\end{tabular}}
\end{table*}

\subsection{Competing algorithms and parameter settings}
To comprehensively evaluate the performance of TRUST-TAEA, several representative large-scale multi-objective evolutionary algorithms (LSMOEAs) are selected as competing methods. These competitors cover different design paradigms, including cooperative co-evolution, reference-direction guidance, population improvement, neural-network training, and sampling strategies. Therefore, they provide a comprehensive basis for assessing the relative performance of the proposed algorithm. The compared algorithms include DGEA \cite{bib38}, FLEA \cite{bib39}, PCPSO \cite{bib40}, LERD \cite{bib41}, NNCSO \cite{bib42}, MSCSO \cite{bib43}, and AMRCSO \cite{bib44}.

To ensure a fair comparison, all algorithms are tested with the same population size and maximum number of generations, and each algorithm is independently executed multiple times under the same experimental environment to reduce the influence of randomness. Specifically, the population size is set to 100, the maximum number of generations is set to 500, and each test instance is independently run 20 times. The final results are compared based on the mean values over multiple independent runs, with standard deviations reported when necessary. For TRUST-TAEA, the basic evolutionary parameters are set to $F=0.5$ and $\textit{CR}=0.9$. The remaining key parameters are determined according to preliminary experiments and further validated through sensitivity analysis in Section~4.4. All algorithms are run under the same function-evaluation budget to avoid comparison bias caused by inconsistent termination conditions.
\begin{figure*}
	\centering
	\includegraphics[width=0.85\linewidth]{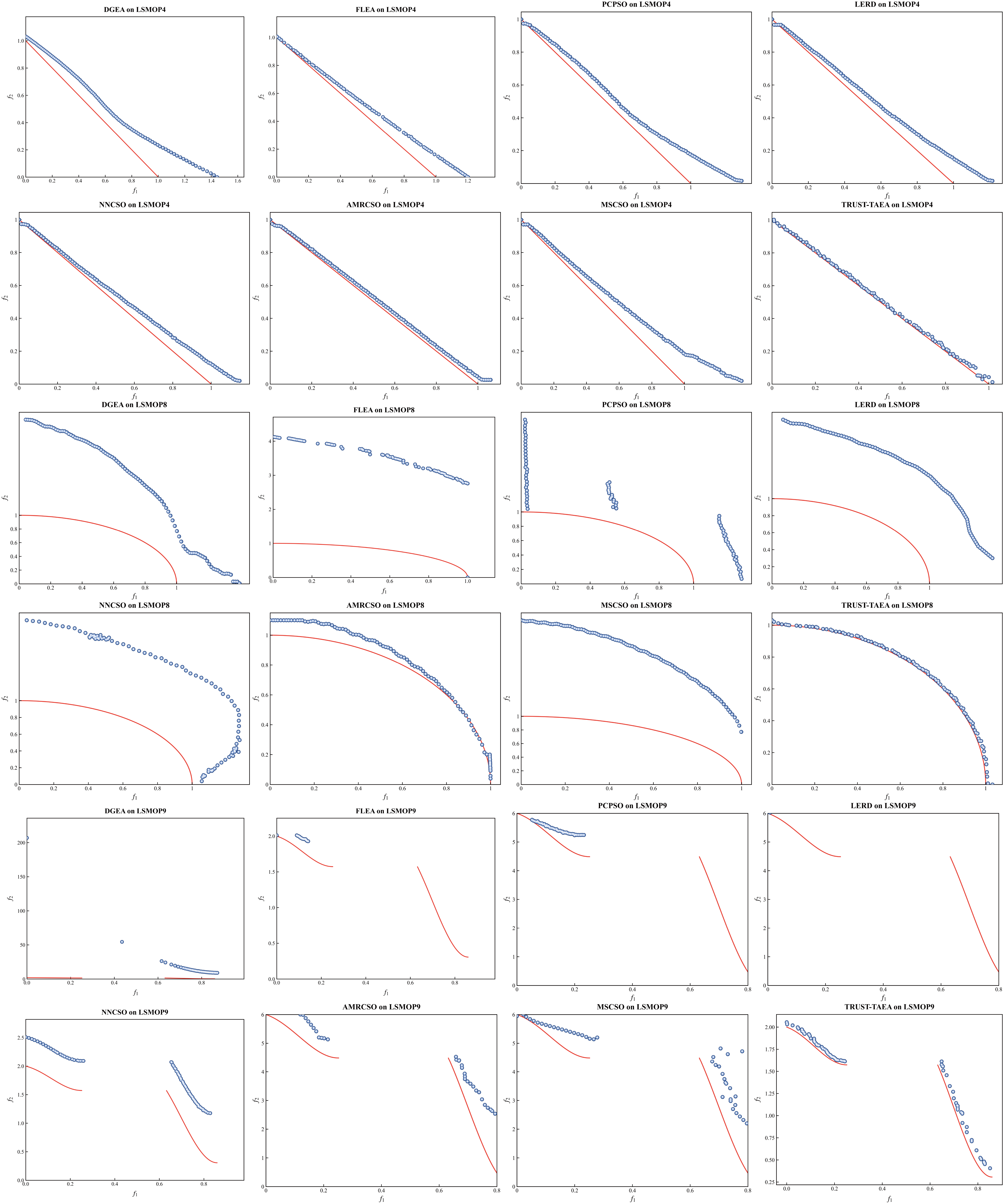} 
	\caption{The true Pareto fronts obtained for LSMOP4, LSMOP8, and LSMOP9 with 5000 decision variables under the 2-objective setting.}
	\label{fig5}
\end{figure*}
\begin{figure*}
	\centering
	\includegraphics[width=0.9\linewidth]{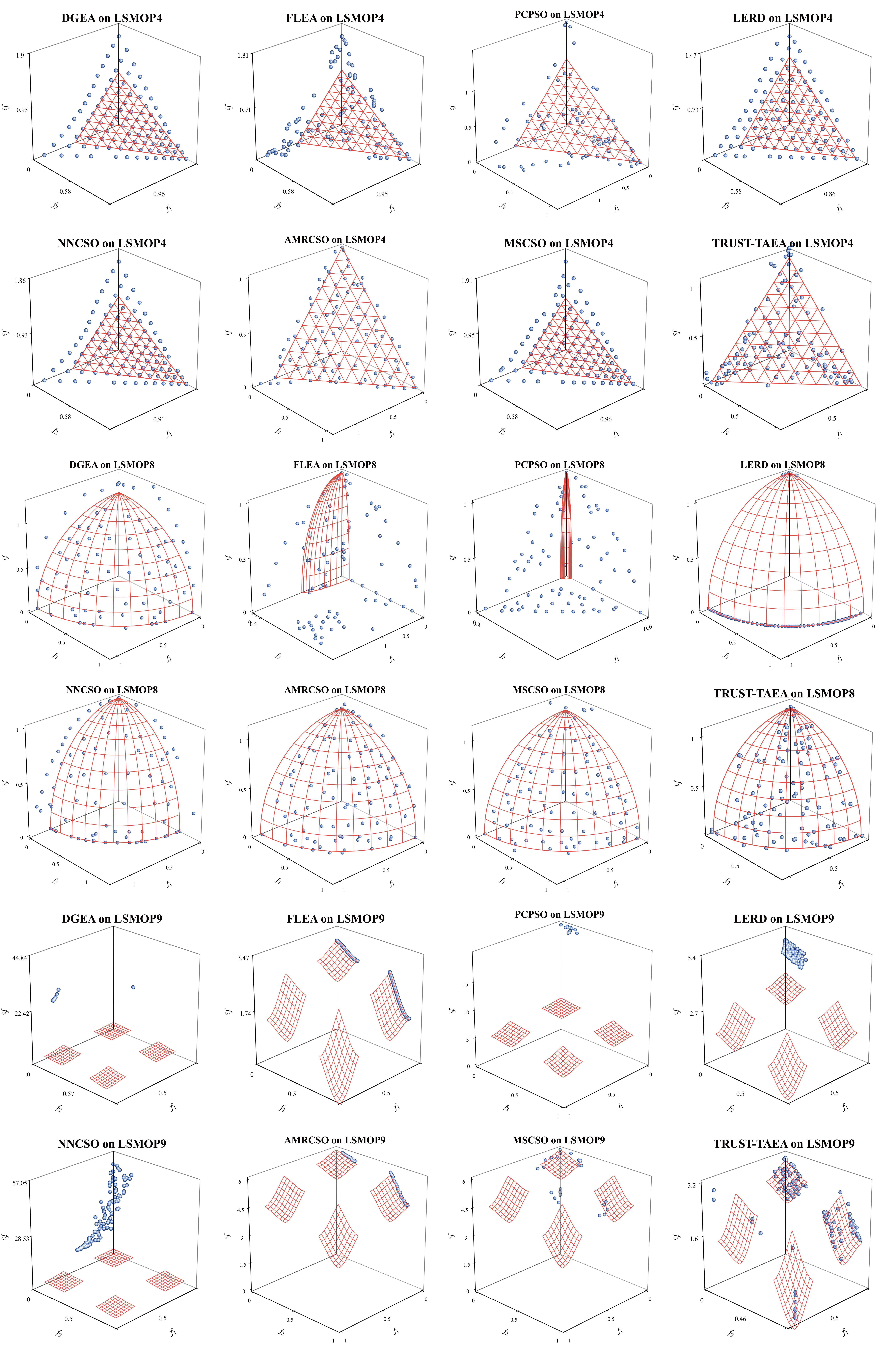} 
	\caption{The true Pareto fronts obtained for LSMOP4, LSMOP8, and LSMOP9 with 5000 decision variables under the 3-objective setting.}
	\label{fig6}
\end{figure*}

\subsection{Performance comparison with SLMOEAs}
Table~\ref{Tab1} shows that TRUST-TAEA achieves the best overall HV performance on the LSMOP test suite. According to the statistical summary in the last row, it significantly outperforms DGEA, FLEA, PCPSO, LERD, NNCSO, AMRCSO, and MSCSO in 65, 66, 66, 67, 67, 65, and 66 out of 72 cases, respectively, while being inferior in only 7, 6, 6, 4, 4, 6, and 5 cases. In the two-objective setting, TRUST-TAEA consistently attains the best HV values on most problems and remains stable as the decision dimension increases, with particularly clear advantages on LSMOP1, LSMOP2, LSMOP4, LSMOP5, LSMOP6, and LSMOP8. It is also notable that on difficult instances such as LSMOP3 and LSMOP7, where many competing algorithms yield zero HV values, TRUST-TAEA still maintains positive dominated volume. In the three-objective setting, its advantages are especially evident on LSMOP1, LSMOP3, LSMOP5, LSMOP6, and LSMOP7, where it often achieves substantially higher HV values than the competing algorithms, indicating stronger front coverage and robustness under challenging high-dimensional conditions. Nevertheless, TRUST-TAEA is less competitive on some three-objective instances, particularly LSMOP2 and LSMOP4, where several competing algorithms obtain larger HV values. Overall, the aforementioned results indicate that TRUST-TAEA performs with exceptional efficiency in addressing highly challenging, large-scale problems prone to Pareto front degradation or insufficient coverage.

Table~\ref{Tab2} shows that TRUST-TAEA achieves the best overall IGD$^+$ performance on the LSMOP test suite. According to the statistical summary in the last row, it significantly outperforms DGEA, FLEA, PCPSO, LERD, NNCSO, AMRCSO, and MSCSO in 65, 72, 71, 68, 62, 64, and 61 out of 72 cases, respectively, while being inferior in only 7, 0, 1, 4, 10, 8, and 9 cases. In the two-objective setting, TRUST-TAEA generally yields the smallest IGD$^+$ values on most problems and remains stable as the decision dimension increases, with especially clear advantages on LSMOP1, LSMOP3, LSMOP5, LSMOP6, LSMOP7, and LSMOP8. In the three-objective setting, its superiority becomes more evident on difficult instances such as LSMOP1, LSMOP3, LSMOP5, LSMOP6, LSMOP7, and LSMOP8, where many competing algorithms exhibit much larger errors or unstable results, whereas TRUST-TAEA consistently maintains relatively small IGD$^+$ values. Nevertheless, TRUST-TAEA is not always the best on every problem; for example, some competing algorithms achieve smaller IGD$^+$ values on parts of LSMOP2 and LSMOP4. Overall, the IGD$^+$ results indicate that TRUST-TAEA has strong front-approximation ability and good robustness, particularly on large-scale instances where competing algorithms tend to suffer from severe convergence loss.
\begin{figure*}
	\centering
	\includegraphics[width=0.8\linewidth]{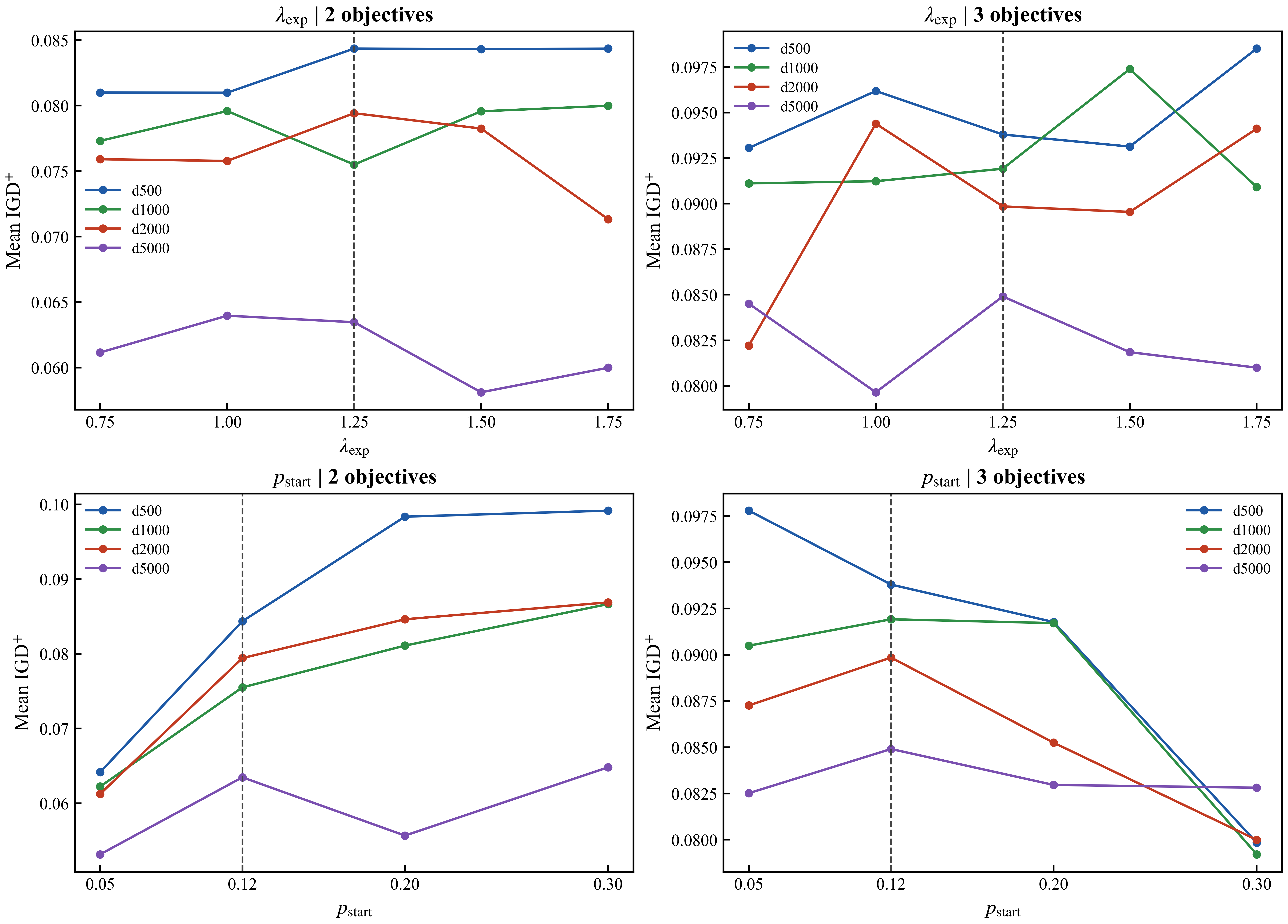} 
	\caption{Overall sensitivity of TRUST-TAEA to $ \lambda_\mathrm{exp} $ and $ p_\mathrm{start} $.}
	\label{fig7}
\end{figure*}
\begin{figure*}
	\centering
	\includegraphics[width=0.8\linewidth]{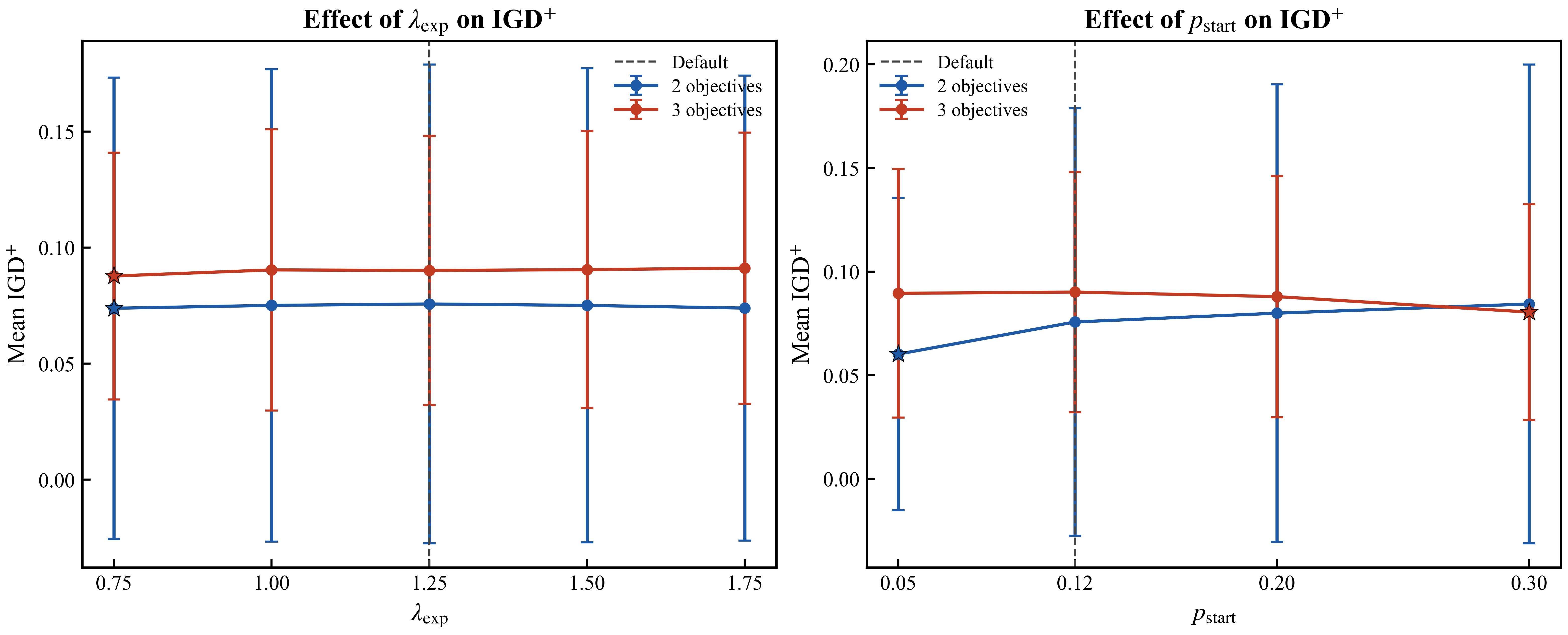} 
	\caption{Dimension-wise sensitivity of TRUST-TAEA to $ \lambda_\mathrm{exp} $ and $ p_\mathrm{start} $.}
	\label{fig8}
\end{figure*}

The above quantitative results verify the effectiveness of TRUST-TAEA at the metric level, but numerical indicators alone cannot fully reflect the differences among algorithms in terms of front-approximation patterns. Therefore, Pareto front visualizations are further provided below. Fig~\ref{fig5} and ~\ref{fig6} illustrate the spatial distribution characteristics of the final approximate solution sets obtained by all algorithms on representative 5000-dimensional LSMOP test problems (LSMOP4, LSMOP8, and LSMOP9) under both the two-objective and three-objective settings. Specifically, AMRCSO shows strong competitiveness in the two-objective case, while in the three-objective case, both AMRCSO and MSCSO are able to obtain high-quality distributed solution sets. However, in comparison, the solution sets produced by TRUST-TAEA exhibit a better overall balance between convergence and diversity. These visual results further verify the effectiveness of the archive-trustworthiness-guided grouping search and anchor-based compensatory mechanism, indicating that these strategies can simultaneously improve the accuracy of approximating the true Pareto front and the uniformity of solution distribution when handling LSMOPs.

\subsection{Parameter sensitivity analysis}
Sensitivity analysis was conducted for $\lambda_{\mathrm{exp}}$ and $p_{\mathrm{start}}$. As shown in Figs.~\ref{fig7} and \ref{fig8}, $\lambda_{\mathrm{exp}}$ has only a mild effect on algorithm performance. For the two-objective problems, the best average IGD$^+$ value is obtained when $\lambda_{\mathrm{exp}}=0.75$, reaching 0.07384, which is approximately $2.4\%$ lower than the default value of 0.07568 obtained at $\lambda_{\mathrm{exp}}=1.25$. For the three-objective problems, the best result is also achieved at $\lambda_{\mathrm{exp}}=0.75$, with an average IGD$^+$ value of 0.08772, corresponding to an improvement of approximately $2.7\%$ over the default value of 0.09012. Moreover, the trends remain generally consistent across different decision dimensions, indicating that TRUST-TAEA is relatively robust to $\lambda_{\mathrm{exp}}$.

In contrast, $p_{\mathrm{start}}$ has a stronger and clearly objective-dependent effect. For the two-objective problems, the best average IGD$^+$ value is obtained when $p_{\mathrm{start}}=0.05$, where it decreases to 0.06019, representing an improvement of approximately $20.5\%$ over the default setting of $p_{\mathrm{start}}=0.12$. For the three-objective problems, the best setting shifts to $p_{\mathrm{start}}=0.30$, yielding an average IGD$^+$ value of 0.08046, which is approximately $10.7\%$ lower than the default value of 0.09012. This result indicates that earlier probe activation is more beneficial for two-objective search, whereas delayed activation is more suitable for three-objective search. Fig.~\ref{fig8} further shows that this trend remains stable across different decision dimensions.

\begin{figure*}
	\centering
	\includegraphics[width=0.8\linewidth]{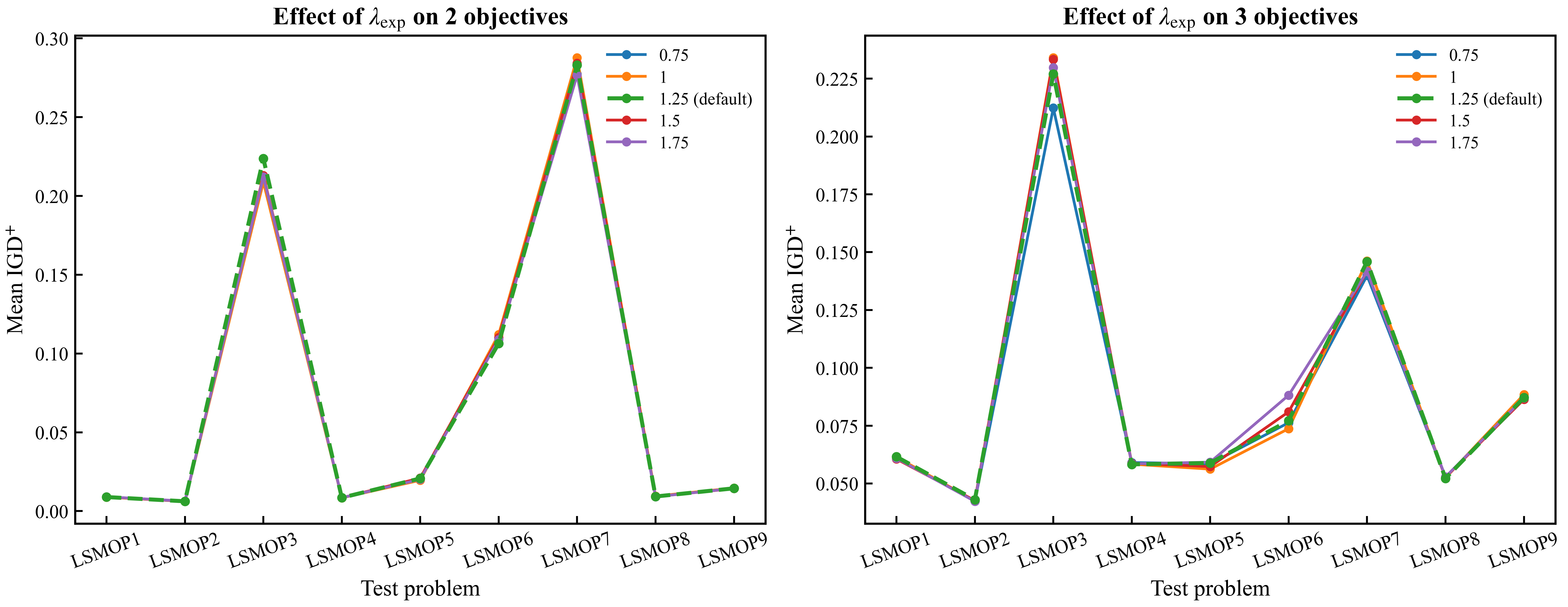} 
	\caption{Problem-wise sensitivity of TRUST-TAEA to $ \lambda_\mathrm{exp} $ on LSMOPs.}
	\label{fig9}
\end{figure*}

\begin{figure*}
	\centering
	\includegraphics[width=0.8\linewidth]{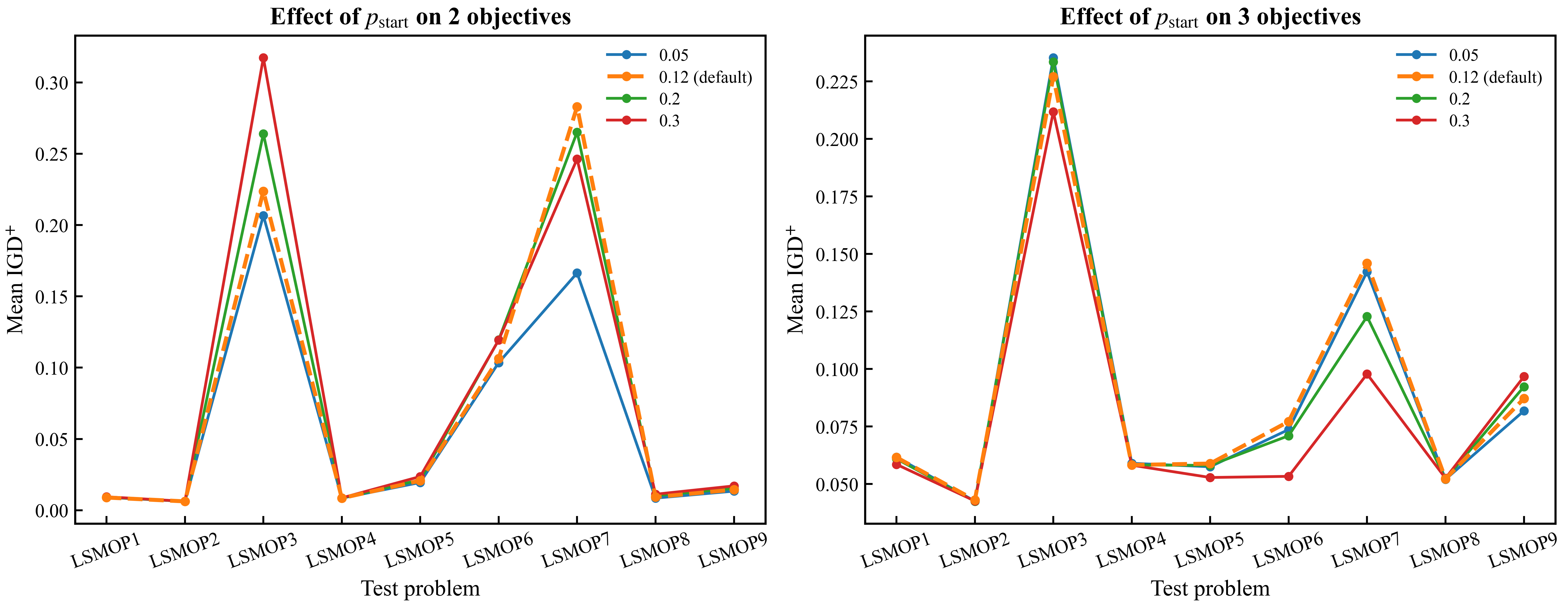} 
	\caption{Problem-wise sensitivity of TRUST-TAEA to $ p_\mathrm{start} $ on LSMOPs.}
	\label{fig10}
\end{figure*}

The problem-wise results in Figs.~\ref{fig9} and \ref{fig10} further support these observations. For $\lambda_{\mathrm{exp}}$, most LSMOP instances show only small performance variations, with relatively stronger sensitivity observed only in a few cases, such as LSMOP3, LSMOP6, and LSMOP7. By contrast, $p_{\mathrm{start}}$ exhibits clearer structural differences. Most two-objective problems prefer smaller values, whereas most three-objective problems prefer larger values, although several functions still show problem-specific preferences. Overall, $\lambda_{\mathrm{exp}}$ acts as a mild tuning parameter, whereas $p_{\mathrm{start}}$ is the key sensitive parameter in TRUST-TAEA. Therefore, under a unified parameter setting, $\lambda_{\mathrm{exp}}$ can be kept at its default value, while $p_{\mathrm{start}}$ should preferably be adjusted according to the number of objectives.

\subsection{Application to a practical microgrid scheduling problem}
To further verify the applicability of the proposed algorithm beyond the LSMOP benchmark suite, a day-ahead scheduling problem of a grid-connected microgrid is considered as the practical application case. Following the general modeling idea in \cite{bib45}, the system includes dispatchable generation, battery storage, demand response, and renewable energy. Unlike \cite{bib45}, which mainly considers cost--emission optimization under uncertainty with clustering-based demand-side management, this study adopts a deterministic representative-day scenario and reformulates the problem as a large-scale three-objective optimization model.

Let $t=1,\ldots,T$ denote the scheduling periods and $\Delta t = 24/T$ be the duration of each period. The decision vector is defined as
\begin{flalign}\label{50}
	& \boldsymbol{x} =
	\left[
	P_1^g, \ldots ,P_T^g,\,
	P_1^{\textit{ch}}, \ldots ,P_T^{\textit{ch}},\,
	P_1^{\textit{dis}}, \ldots ,P_T^{\textit{dis}},\,
	\right. &\notag \\
	&\left.
	\qquad
	P_1^{\textit{dr}}, \ldots ,P_T^{\textit{dr}},\,
	P_1^{\textit{cur}}, \ldots ,P_T^{\textit{cur}}
	\right] \in \mathbb{R}^{5T}, &
\end{flalign}
where $P_t^g$, $P_t^{\textit{ch}}$, $P_t^{\textit{dis}}$, $P_t^{\textit{dr}}$, and $P_t^{\textit{cur}}$ denote the generator output, battery charging power, battery discharging power, demand response amount, and renewable curtailment, respectively. Under the default setting $T=96$, the decision dimension is $D=5T=480$.

The grid purchase power is determined by the following power-balance equation:
\begin{flalign}\label{51}
	& P_t^{\textit{grid}} =
	L_t - P_t^{\textit{dr}} + P_t^{\textit{ch}} + P_t^{\textit{cur}}
	- \left(P_t^g + P_t^{\textit{dis}} + P_t^{\textit{re}}\right), & \notag \\
	& P_t^{\textit{grid}} \ge 0, &
\end{flalign}
where $L_t$ and $P_t^{\textit{re}}$ are the load demand and available renewable power at period $t$, respectively.

The three objective functions are formulated as
\begin{flalign}\label{52,53,54}
	&\min \; f_1(\boldsymbol {x}) = \Delta t \sum_{t=1}^{T} \left( c_t^{\textit{grid}} P_t^{\textit{grid}} + c^{g} P_t^{g} + c^{b}\left(P_t^{\textit{ch}}+P_t^{\textit{dis}}\right) \right. & \notag \\ &\qquad\qquad\qquad\;\;\;\left. + c^{\textit{dr}} P_t^{\textit{dr}} + c^{\textit{cur}} P_t^{\textit{cur}}\right),& \\ &\min \; f_2(\boldsymbol {x})=\Delta t \sum_{t=1}^{T}\left(e_t^{\textit{grid}} P_t^{\textit{grid}}+e^{g} P_t^{g}\right),& \\ &\min \; f_3(\boldsymbol {x})=\frac{1}{T-1}\sum_{t=2}^{T}\left|P_t^{\textit{grid}}-P_{t-1}^{\textit{grid}}\right|.& 
\end{flalign}

The model is subject to the following operational constraints:
\begin{flalign}\label{55,56,57,58,59,60,61}
	& 0 \le P_t^g \le P^{g}_{\max}, \;
	0 \le P_t^{\textit{ch}} \le P^{\textit{ch}}_{\max}, \;
	0 \le P_t^{\textit{dis}} \le P^{\textit{dis}}_{\max},&
	\\
	&0 \le P_t^{\textit{dr}} \le P_t^{\textit{dr},\max}, \;
	0 \le P_t^{\textit{cur}} \le P_t^{\textit{re}},&
	\\
	&|P_t^g - P_{t-1}^g| \le R^g, \; t=2,\ldots,T,&
	\\
	& E_t= E_{t-1} + \eta^{\textit{ch}}\Delta t\, P_t^{\textit{ch}}
	- \frac{\Delta t}{\eta^{\textit{dis}}} P_t^{\textit{dis}},&
	\\
	& E_{\min} \le E_t \le E_{\max},\; t=1,\ldots,T,&
	\\
	& \Delta t \sum_{t=1}^{T} P_t^{\textit{dr}} \le E^{\textit{dr}}_{\max},&
	\\
	& |E_T - E_0| \le \varepsilon.&
\end{flalign}

Compared with \cite{bib45}, the improved model used in this study has three main differences. First, a grid-fluctuation objective is introduced, so that the dispatch scheme simultaneously accounts for economic cost, environmental impact, and operational smoothness. Second, the problem is represented directly by a high-dimensional continuous decision vector, making it more suitable for evaluating LSMOEAs. Third, renewable curtailment and terminal SOC balance are explicitly incorporated. Therefore, this application problem preserves the main engineering characteristics of microgrid day-ahead scheduling while providing a more challenging practical test for the proposed algorithm.

\begin{table}[htbp]
	\centering
	\caption{Average IGD$^+$ results of the compared algorithms on the three-objective microgrid dispatch problem over 10 independent runs.}
	\label{Tab3}
	\begin{tabular*}{\columnwidth}{@{\extracolsep{\fill}}lcc@{}}
		\toprule
		Algorithm & IGD$^+$ (Mean $\pm$ Std.) & Rank \\
		\midrule
		DGEA       & 1.1305 $\pm$ 0.0747 & 7 \\
		FLEA       & 0.4700 $\pm$ 0.0246 & 3 \\
		PCPSO      & 1.1695 $\pm$ 0.0714 & 8 \\
		LERD       & 0.3112 $\pm$ 0.0469 & 2 \\
		NNCSO      & 0.5837 $\pm$ 0.0550 & 5 \\
		AMRCSO     & 0.4711 $\pm$ 0.0536 & 4 \\
		MSCSO      & 0.7409 $\pm$ 0.0483 & 6 \\
		TRUST-TAEA & \textbf{0.0977 $\pm$ 0.0248} & \textbf{1} \\
		\bottomrule
	\end{tabular*}
\end{table}

\begin{figure}
	\centering
	\includegraphics[width=\linewidth]{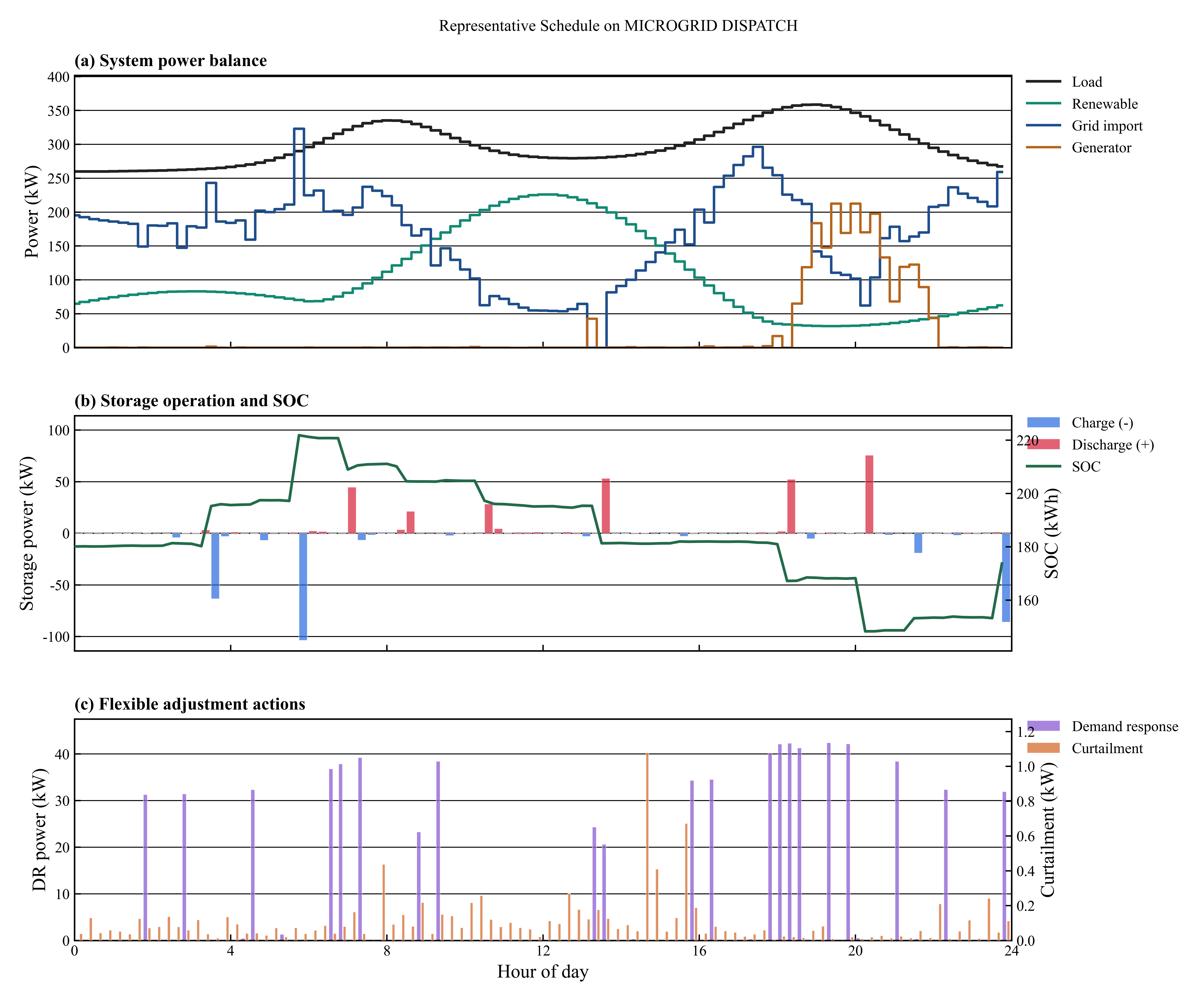}
	\caption{Representative dispatch schedule obtained by TRUST-TAEA in its best run on the three-objective microgrid dispatch problem. In panel (a), the black, green, blue, and brown curves denote the load demand, available renewable power, grid purchase power, and generator output, respectively. In panel (b), the blue and red bars represent battery charging and discharging power, respectively, where charging power is plotted on the negative side to indicate energy absorption, and the green curve on the secondary axis denotes the battery state of charge (SOC). In panel (c), the purple bars denote the demand response amount and the orange bars denote renewable curtailment; curtailment is plotted on the secondary axis because its magnitude is much smaller than that of demand response.}
	\label{fig11}
\end{figure}

As shown in Table~\ref{Tab3}, TRUST-TAEA achieves the best average IGD$^+$ value of 0.0977, whereas the second-best algorithm, LERD, obtains 0.3112. This corresponds to a 68.60\% reduction in average IGD$^+$ relative to LERD. The advantage is also evident when compared with FLEA and AMRCSO, whose average IGD$^+$ values are 0.4700 and 0.4711, respectively. By contrast, DGEA and PCPSO both produce average IGD$^+$ values above 1.13, indicating a much larger distance from the empirical reference front. These results show that TRUST-TAEA achieves the best overall convergence performance on this practical scheduling problem.

To further explain the obtained dispatch behavior, Fig.~\ref{fig11} presents a representative schedule from the best run of TRUST-TAEA. Panel (a) shows the system-level power balance, panel (b) shows the battery charging/discharging actions and SOC trajectory, and panel (c) reports the flexible adjustment actions, including demand response and renewable curtailment. In this run, the selected solution reaches $(f_1,f_2,f_3)=(3016.83,\,2841.26,\,20.25)$ with zero penalty, indicating full feasibility.

Over the 24-hour horizon, the total load demand and available renewable energy are 7181.71~kWh and 2461.64~k\-Wh, respectively. The obtained schedule uses 4003.94~kWh of grid electricity and 534.48~kWh of generator output, while coordinating 185.89~kWh of demand response and 79.82/78.03~kWh of battery charging/discharging. Meanwhile, renewable curtailment is only 2.45~kWh, accounting for approximately 0.10\% of the available renewable energy.

It can also be observed from Fig.~\ref{fig11} that the generator is mainly activated during the evening peak-load period, demand response is concentrated in high-demand intervals, and the battery performs targeted energy shifting rather than excessive cycling. The SOC varies between 148.33~kWh and 221.81~kWh and returns to 173.70~kWh at the end of the scheduling horizon. This value remains close to the initial SOC of 180~kWh and satisfies the terminal energy constraint. These results indicate that TRUST-TAEA not only achieves the best convergence performance among the compared algorithms, but also produces a feasible and well-structured dispatch strategy for the practical microgrid scheduling problem.

\section{Conclusion}\label{sec5}
This paper proposed TRUST-TAEA, a trustworthiness-guided two-archive evolutionary algorithm for solving LSM\-OPs. The proposed method extends the conventional two-archive framework by introducing an archive trustworthiness indicator to regulate offspring generation and archive feedback during evolution. Based on this indicator, variable-grouping sparse search, anchor-probing compensatory search, and checkpoint-based archive stabilization are developed to improve search efficiency, front coverage, and evolutionary stability in high-dimensional decision spaces. Experimental studies on 2, 3-objective LSMOP benchmark problems with decision dimensions ranging from 500 to 5000 show that TRUST-TAEA achieves competitive, and in most cases superior, performance compared with the selected algorithms in terms of HV and IGD$^+$. These results indicate that the proposed trustworthiness-guided control mechanism can effectively coordinate convergence-oriented search and diversity maintenance under limited evaluation budgets. In addition, the parameter sensitivity analysis shows that the main parameters of TRUST-TAEA are reasonably robust within appropriate ranges. To further examine its practical applicability, TRUST-TAEA was also applied to a three-objective day-ahead scheduling problem of a grid-connected microgrid. The results show that the proposed algorithm achieves the best IGD$^+$ value among the compared algorithms and produces a feasible dispatch schedule that balances operation cost, carbon emissions, and grid-power fluctuation. The obtained schedule coordinates grid purchase, dispatchable generation, battery charging and discharging, demand response, and renewable curtailment in a structured manner. This case study demonstrates that TRUST-TAEA is not only effective on benchmark problems, but also has practical potential for solving high-dimensional engineering optimization problems with multiple conflicting objectives.

In summary, this work provides an effective approach for enhancing two-archive evolutionary optimization in large-scale multi-objective problems. Future work will focus on extending the proposed framework to more complex scenarios, such as constrained LSMOPs, dynamic LSMOPs, uncertain environments, and adaptive parameter control.

\section*{Credit authorship contribution statement}
\textbf{Junyi Cui:} Writing – original draft, Software, Methodology, Data curation, Formal analysis. \textbf{Chao Min}: Formal analysis, Investigation, Resources, Writing – review \& editing, Validation. \textbf{Stanis{\l}aw Mig{\'o}rski}: Validation, Formal analysis, Visualization, Supervision. \textbf{Binrong Wang}: Formal analysis, Validation, Supervision. \textbf{Yonglan Xie}: Validation, Writing – review \& editing, Supervision.

\section*{Declaration of competing interest}
The authors declare that they have no known competing financial interests or personal relationships that could have appeared to influence the work reported in this paper.

\section*{Acknowledgments}
This work was supported by Sichuan Provincial Science and Technology Project (No.2025NSFTD0016) and the Open Bidding for Selecting the Best Candidates Project of Southwest Petroleum University (2024CXJB11).

\printcredits

\bibliographystyle{elsarticle-num}

\bibliography{cas-refs}



\end{document}